\pdfoutput=1
\documentclass[letterpaper,11pt]{article}
\usepackage[margin=1in]{geometry}
\usepackage[CJKbookmarks=true,
            bookmarksnumbered=true,
            bookmarksopen=true,
            colorlinks=true,
            citecolor=red,
            linkcolor=blue,
            anchorcolor=red,
            urlcolor=blue
            ]{hyperref}
\usepackage[title]{appendix}
\usepackage{graphicx}
\usepackage{booktabs}
\usepackage{algorithm}
\usepackage{algorithmic}
\usepackage{amsmath,amsthm,amsfonts,amssymb}
\usepackage{color}
\usepackage{enumitem}
\usepackage{bm}
\usepackage{multirow}
\usepackage{bbm}
\usepackage{subfiles}
\usepackage{xspace}
\usepackage{caption}
\usepackage{mathtools} 
\usepackage{tikz}
\usepackage{subfig}
 \usepackage{natbib}
 \usepackage{xcolor}

\usetikzlibrary{arrows.meta,patterns,calc}

\newtheorem{theorem}{Theorem}
\newtheorem*{theorem*}{Theorem}
\newtheorem{lemma}{Lemma}

\newtheorem{proposition}{Proposition}

\theoremstyle{definition}
\newtheorem{remark}{Remark}
\newtheorem{definition}{Definition}

\tikzset{global scale/.style={
    scale=#1,
    every node/.append style={scale=#1}
  }
}

\newcommand{\maC}{\mathcal{C}}

\newcommand{\maE}{\mathcal{E}}

\newcommand{\maG}{\mathcal{G}}
\newcommand{\maH}{\mathcal{H}}
\newcommand{\maI}{\mathcal{I}}

\newcommand{\maP}{\mathcal{P}}
\newcommand{\maQ}{\mathcal{Q}}

\newcommand{\maT}{\mathcal{T}}

\newcommand{\prob}[1]{ \mathbb{P}\left[ #1 \right] }
\newcommand{\expect}[1]{ \mathbb{E}\left[ #1 \right] }

\newcommand{\diff}{{\rm d}}
\newcommand{\Expect}{\mathbb{E}}

\newcommand{\vn}{\mathsf{{v}}}

\newcommand{\pth}[1]{\left( #1 \right)}
\newcommand{\qth}[1]{\left[ #1 \right]}
\newcommand{\sth}[1]{\left\{ #1 \right\}}

\newcommand{\calI}{{\mathcal{I}}}

\newcommand{\ti}{\tilde}
\newcommand{\TV}{\mathsf{TV}}
\newcommand{\Bin}{\mathrm{Bin}}
\newcommand{\Bernoulli}{\mathrm{Bernoulli}}

\newcommand{\sfC}{{\mathsf{C}}}

\newcommand{\sfE}{{\mathsf{E}}}

\newcommand{\sfP}{{\mathsf{P}}}

\newcommand{\sfS}{{\mathsf{S}}}
\newcommand{\sfT}{{\mathsf{T}}}

\newcommand{\sfV}{{\mathsf{V}}}

\newcommand{\HG}{\mathrm{HG}}

\newcommand{\indc}[1]{{\mathbf{1}_{\left\{{#1}\right\}}}}

\newcommand{\ER}{Erd\H{o}s-R\'enyi }

\newcommand{\Aut}{\mathsf{Aut}}

\newcommand{\inj}{\mathsf{inj}}
\newcommand{\vp}{\varphi}

\def\var{\mathrm{Var}}
\def\E{\mathbb{E}}
\def\P{\mathbb{P}}

\usepackage{xcolor}
\definecolor{red1}{RGB}{255, 0, 0}  
\definecolor{red2}{RGB}{220, 47, 2}   
\definecolor{red3}{RGB}{189, 58, 85}
\usetikzlibrary{intersections,calc,decorations.pathreplacing,through,arrows}
\usetikzlibrary{ decorations.markings,positioning}

\title{Sample Complexity of Correlation Detection in \ER graphs}
\author{Dong Huang and Pengkun Yang\thanks{
D.\ Huang and P.\ Yang are with the Department of Statistics and Data Science, Tsinghua
University. P. Yang is supported in part by National Key R\&D Program
of China 2024YFA1015800, Tsinghua
University Dushi Program.
}}
\date{}

\begin{document}

\maketitle

\begin{abstract}
In this paper, we study correlation detection problem between a pair of \ER graphs. Under the null hypothesis, the two graphs are independent, whereas under the alternative, they are edge-correlated through a latent permutation. We consider the partial-observation setting in which only an induced subgraph of each graph is sampled, and determine the sharp information-theoretic threshold for the sample size, including its leading constant. We also establish a general conditional second moment principle and identify a continuous transition in the sharp threshold constant across two different regimes.
\end{abstract}

\tableofcontents

\section{Introduction}\label{sec:intro}

Graph alignment has recently attracted substantial interest. Given two networks defined on a common latent set of entities but observed under an unknown vertex permutation, the goal is to recover the underlying vertex correspondence from their structural similarities. Applications include matching user accounts across social networks~\citep{narayanan2008robust,narayanan2009anonymizing}, homologous proteins across species~\citep{singh2008global,vogelstein2015fast}, corresponding parts of images or 3-D shapes~\citep{berg2005shape,mateus2008articulated}, and related linguistic entities across sentences or languages~\citep{hughes2007lexical}.

However, finding the optimal vertex correspondence between any two graphs needs to solve the NP-hard quadratic assignment problem (QAP), which leads to the worst-case intractability~\citep{pardalos1994quadratic}. Consequently, many recent works focus on the average-case analysis of graph alignment under meaningful random graph models, such as correlated \ER model~\citep{pedarsani2011privacy} and correlated Gaussian Wigner model~\citep{ding2021efficient}.
In this paper, we consider the following statistical model where two random graphs are correlated through a hidden vertex correspondence. Let $\maG(n,p)$ denote the distribution of \ER random graphs on $n$ vertices, where each potential edge is present independently with probability $p\in (0,1)$. We first introduce the following correlated \ER graph model.

\begin{definition}[Correlated \ER graph]
    Let $\mathbf{G}_1,\mathbf{G}_2$ be two random graphs with vertex sets $V(\mathbf{G}_1),V(\mathbf{G}_2)$ such that $|V(\mathbf{G}_1)|  = |V(\mathbf{G}_2)| = n$.
    Let $\pi^*$ denote the latent uniform bijective mapping from $V(\mathbf{G}_1)$ to $V(\mathbf{G}_2)$. We say a pair of graphs $(\mathbf{G}_1,\mathbf{G}_2)$ follows correlated \ER graphs distribution $\maG(n,p,\rho)$ if both marginal distributions are $\maG(n,p)$ and each pair of edges $(uv,\pi^*(u)\pi^*(v))$ for $u,v\in V(\mathbf{G}_1)$ follows the correlated bivariate Bernoulli distribution with correlation coefficient $\rho$.
\end{definition}

Despite the substantial progress on the graph alignment problem, much less attention has been paid to the arguably more basic task of detecting latent correlation between two networks, as opposed to testing whether the graphs are generated independently.  
Following the hypothesis testing framework in~\cite{barak2019nearly}, for two graphs $\mathbf{G}_1,\mathbf{G}_2$, we consider the following problem:
\begin{itemize}
    \item Under the null hypothesis $\maH_0$, $\mathbf{G}_1$ and $\mathbf{G_2}$ are independently generated from the \ER graph model $\maG(n,p)$;
    \item Under the alternative hypothesis $\maH_1$, $\mathbf{G}_1$ and $\mathbf{G}_2$ are generated from the correlated \ER graph model $\maG(n,p,\rho)$;
\end{itemize}
We note that the marginal distributions of $\mathbf{G}_1$ and $\mathbf{G}_2$ are both $\maG(n,p)$ under $\maH_0$ and $\maH_1$; the two hypotheses differ only in the latent vertex permutation $\pi^*$ under $\maH_1$. Consequently, any powerful test must exploit the joint structural information across the two graphs rather than statistics of each graph in isolation.

\subsection{Graph sampling}

There have been substantial prior works on detection problems for a pair of \ER graphs, including information-theoretic analyses~\citep{wu2023testing,ding2023detection} and efficient algorithm analyses~\citep{barak2019nearly,mao2024testing,huang2025testing}. However, most existing methods assume full access to all edges in both $\mathbf{G}_1$ and $\mathbf{G}_2$, an assumption that is often unrealistic in practice. In many applications, only sampled subgraphs of the underlying networks are observable, and graph sampling becomes a natural tool for exploring their structure, with examples including:
\begin{itemize}
    \item \emph{Lack of data.} In social network analysis, the complete network is rarely available, and
    one often works with a sampled subset of vertices or edges~\citep{papagelis2011sampling}.
    \item \emph{Testing costs.} In protein–protein interaction studies, experimentally testing all
    potential interactions is expensive, so sampling-based designs are routinely employed~\citep{stumpf2005subnets}.
    \item \emph{Surveying hidden populations.} In sociological studies, communities within networks are
    frequently explored via sampled ego-networks using techniques such as snowball sampling
    and respondent-driven sampling~\citep{goodman1961snowball,salganik20045,lancichinetti2009community}.
    \item \emph{Visualization.} For very large graphs, sampling a subgraph often provides a more digestible
    representation for visualization and exploratory analysis~\citep{wu2016evaluation}.
\end{itemize}

Graph sampling refers to observing only a partial view of a large network in order to approximate its structural or statistical properties \citep{leskovec2006sampling}. 
Common designs include (i) \emph{vertex sampling}, which selects a subset of vertices and observes the induced subgraph on them, 
(ii) \emph{edge sampling}, which selects a subset of edges and keeps their incident endpoints, and 
(iii) \emph{exploration-based sampling} such as breadth-first search sampling \citep{leskovec2006sampling}, snowball sampling \citep{goodman1961snowball}, or Metropolis--Hastings random-walk sampling \citep{gjoka2010walking}.

In this paper, we consider the induced subgraphs sampling methods for hypothesis testing between $\maH_0$ and $\maH_1$. Specifically, for $\mathbf{G}_1,\mathbf{G}_2$ with vertex sets $V(\mathbf{G}_1),V(\mathbf{G}_2)$ such that $\vn(\mathbf{G}_1) = \vn(\mathbf{G}_2) = n$, where $\vn(G)\triangleq |V(G)|$, we randomly sample two induced subgraphs $G_1,G_2$ with $s$ vertices. An induced subgraph of a graph is formed from a subset of the vertices of the graph, along with all the edges between them from the original graph. Specifically, we first independently select vertex sets $V(G_1)\subseteq V(\mathbf{G}_1)$ and $V(G_2)\subseteq V(\mathbf{G}_2)$ with $\vn(G_1) = \vn(G_2) = s$, and then keep the edges between $V(G_1)$ and $V(G_2)$ from the original graphs (see Figure~\ref{fig:induced-subgraph} for an illustration). Let $\maQ$ and $\maP$ denote the probability measure of sampled subgraphs $(G_1,G_2)$ under $\maH_0$ and $\maH_1$, respectively. We say that a test statistic $\maT(G_1,G_2)$ with a  threshold $\tau\in \mathbb{R}$ achieves detection, if the sum of type I and type II error converges to 0 as $n\to \infty$,\begin{align}\label{eq:strong-detection}
        \lim_{n\to \infty } \qth{\maQ(\maT(G_1,G_2)\ge \tau)+{\maP}(\maT(G_1,G_2)< \tau)} = 0;
    \end{align}

\begin{figure}[htbp]
\centering
\begin{tikzpicture}[
  x=1cm,y=1cm,
  scale=0.45, transform shape,
  rednode/.style={circle,draw=none,fill=red!85,minimum size=3mm,inner sep=0pt},
  grynode/.style={circle,draw=none,fill=gray!35,minimum size=3mm,inner sep=0pt},
  bedge/.style={draw=black,line width=0.9pt,opacity=0.85},
  gedge/.style={draw=gray!60,line width=0.9pt,opacity=0.35},
  title/.style={font=\huge}
]

\pgfmathsetmacro{\ys}{0.65}

\begin{scope}[shift={(0,0)}]
  \node[title] at (0,4.2) {Original Graph $\mathbf{G}_1$ ($n=10$)};

  \node[rednode] (g1-1)  at ( 3.80, 0.00) {};
  \node[rednode] (g1-2)  at ( 3.07, { 2.24*\ys}) {};
  \node[rednode] (g1-3)  at ( 1.17, { 3.61*\ys}) {};
  \node[rednode] (g1-4)  at (-1.17, { 3.61*\ys}) {};
  \node[rednode] (g1-5)  at (-3.07, { 2.24*\ys}) {};
  \node[rednode] (g1-6)  at (-3.80, 0.00) {};
  \node[rednode] (g1-7)  at (-3.07, {-2.24*\ys}) {};

  \node[grynode] (g1-8)  at (-1.17, {-3.61*\ys}) {};
  \node[grynode] (g1-9)  at ( 1.17, {-3.61*\ys}) {};
  \node[grynode] (g1-10) at ( 3.07, {-2.24*\ys}) {};

  \draw[bedge] (g1-4)--(g1-5);
  \draw[bedge] (g1-5)--(g1-6);
  \draw[bedge] (g1-6)--(g1-7);
  \draw[bedge] (g1-2)--(g1-4);
  \draw[bedge] (g1-3)--(g1-4);
  \draw[bedge] (g1-3)--(g1-5);
  \draw[bedge] (g1-1)--(g1-4);
  \draw[bedge] (g1-1)--(g1-7);
  \draw[bedge] (g1-2)--(g1-6);
  \draw[bedge] (g1-4)--(g1-7);

  \draw[gedge] (g1-4)--(g1-8);
  \draw[gedge] (g1-5)--(g1-8);
  \draw[gedge] (g1-2)--(g1-8);

  \draw[gedge] (g1-4)--(g1-9);
  \draw[gedge] (g1-5)--(g1-9);
  \draw[gedge] (g1-3)--(g1-9);
  \draw[gedge] (g1-6)--(g1-9);

  \draw[gedge] (g1-5)--(g1-10);
  \draw[gedge] (g1-7)--(g1-10);
  \draw[gedge] (g1-9)--(g1-10);
  \draw[gedge] (g1-8)--(g1-9);
\end{scope}

\begin{scope}[shift={(10.2,0)}]
  \node[title] at (0,4.2) {Original Graph $\mathbf{G}_2$ ($n=10$)};

  \node[rednode] (g2-1)  at ( 3.80, 0.00) {};
  \node[rednode] (g2-2)  at ( 3.07, { 2.24*\ys}) {};
  \node[rednode] (g2-3)  at ( 1.17, { 3.61*\ys}) {};
  \node[rednode] (g2-4)  at (-1.17, { 3.61*\ys}) {};
  \node[rednode] (g2-5)  at (-3.07, { 2.24*\ys}) {};
  \node[rednode] (g2-6)  at (-3.80, 0.00) {};
  \node[rednode] (g2-7)  at (-3.07, {-2.24*\ys}) {};

  \node[grynode] (g2-8)  at (-1.17, {-3.61*\ys}) {};
  \node[grynode] (g2-9)  at ( 1.17, {-3.61*\ys}) {};
  \node[grynode] (g2-10) at ( 3.07, {-2.24*\ys}) {};

  \draw[bedge] (g2-1)--(g2-2);
  \draw[bedge] (g2-1)--(g2-3);
  \draw[bedge] (g2-1)--(g2-6);
  \draw[bedge] (g2-2)--(g2-5);
  \draw[bedge] (g2-2)--(g2-7);
  \draw[bedge] (g2-3)--(g2-7);
  \draw[bedge] (g2-4)--(g2-7);
  \draw[bedge] (g2-6)--(g2-7);

  \draw[gedge] (g2-3)--(g2-8);
  \draw[gedge] (g2-2)--(g2-8);
  \draw[gedge] (g2-1)--(g2-8);

  \draw[gedge] (g2-8)--(g2-9);
  \draw[gedge] (g2-9)--(g2-10);
  \draw[gedge] (g2-8)--(g2-10);

  \draw[gedge] (g2-3)--(g2-9);
  \draw[gedge] (g2-3)--(g2-10);
  \draw[gedge] (g2-7)--(g2-10);
  \draw[gedge] (g2-1)--(g2-10);
\end{scope}

\begin{scope}[shift={(0,-8.2)}]
  \node[title] at (0,4.2) {Induced Subgraph $G_1$ ($s=7$)};

  \node[rednode] (g1s-1) at ( 3.80, 0.00) {};
  \node[rednode] (g1s-2) at ( 3.07, { 2.24*\ys}) {};
  \node[rednode] (g1s-3) at ( 1.17, { 3.61*\ys}) {};
  \node[rednode] (g1s-4) at (-1.17, { 3.61*\ys}) {};
  \node[rednode] (g1s-5) at (-3.07, { 2.24*\ys}) {};
  \node[rednode] (g1s-6) at (-3.80, 0.00) {};
  \node[rednode] (g1s-7) at (-3.07, {-2.24*\ys}) {};

  \draw[bedge] (g1s-4)--(g1s-5);
  \draw[bedge] (g1s-5)--(g1s-6);
  \draw[bedge] (g1s-6)--(g1s-7);
  \draw[bedge] (g1s-2)--(g1s-4);
  \draw[bedge] (g1s-3)--(g1s-4);
  \draw[bedge] (g1s-3)--(g1s-5);
  \draw[bedge] (g1s-1)--(g1s-4);
  \draw[bedge] (g1s-1)--(g1s-7);
  \draw[bedge] (g1s-2)--(g1s-6);
  \draw[bedge] (g1s-4)--(g1s-7);
\end{scope}

\begin{scope}[shift={(10.2,-8.2)}]
  \node[title] at (0,4.2) {Induced Subgraph $G_2$ ($s=7$)};

  \node[rednode] (g2s-1) at ( 3.80, 0.00) {};
  \node[rednode] (g2s-2) at ( 3.07, { 2.24*\ys}) {};
  \node[rednode] (g2s-3) at ( 1.17, { 3.61*\ys}) {};
  \node[rednode] (g2s-4) at (-1.17, { 3.61*\ys}) {};
  \node[rednode] (g2s-5) at (-3.07, { 2.24*\ys}) {};
  \node[rednode] (g2s-6) at (-3.80, 0.00) {};
  \node[rednode] (g2s-7) at (-3.07, {-2.24*\ys}) {};

  \draw[bedge] (g2s-1)--(g2s-2);
  \draw[bedge] (g2s-1)--(g2s-3);
  \draw[bedge] (g2s-1)--(g2s-6);
  \draw[bedge] (g2s-2)--(g2s-5);
  \draw[bedge] (g2s-2)--(g2s-7);
  \draw[bedge] (g2s-3)--(g2s-7);
  \draw[bedge] (g2s-4)--(g2s-7);
  \draw[bedge] (g2s-6)--(g2s-7);
\end{scope}

\end{tikzpicture}
\caption{An example of induced subgraph sampling with $n=10$ and $s=7$.}
\label{fig:induced-subgraph}
\end{figure}

Intuitively, to meet the detection criterion the two sampled subgraphs must have an overlap, which is equivalent to requiring that the sampled vertex sets intersect. If we sample $s$ vertices from a population of $n$ vertices twice without replacement, the expected size of the intersection of the two vertex sets is $s^{2}/n$. Hence, for the overlap to be of constant order, the sample size must satisfy $s \gtrsim \sqrt{n}$, which yields the minimal scaling for successful detection. Indeed, applying a union bound introduces an additional $\log n$ factor in the detection threshold, leading to the requirement $s \gtrsim \sqrt{n \log n}$. See Section~\ref{sec:upbd} for further details.

\subsection{Main results}

In this subsection, we present the main results on information-theoretic thresholds. We first introduce some notations for the presentation of the main theorems. Throughout the paper, we assume that $0<\rho<1, 0<p\le \frac{1}{2}$, and the sample size $s\le n$. We denote the bivariate distribution of a pair of Bernoulli random variables with means $p_1,p_2$, and correlation coefficient $\rho$ as $\mathrm{Bern}(p_1,p_2,\rho)$. In the correlated \ER model, a pair of correlated edges $(uv,\pi^*(u)\pi^*(v))$ follows bivariate Bernoulli distribution $\mathrm{Bern}(p,p,\rho)$, where $u,v\in V(\mathbf{G}_1)$. Specifically, for $\mathrm{Bern}(p,p,\rho)$,
we denote the following probability mass function:\begin{align*}
    p_{11}\triangleq p^2+\rho p (1-p),\quad p_{10} = p_{01} \triangleq (1-\rho) p (1-p),\quad p_{00} \triangleq (1-p)^2+\rho  p (1-p).
\end{align*}
We define $\gamma \triangleq \frac{p^2+\rho p (1-p)}{p^2}-1 = \frac{\rho (1-p)}{p}$, which quantifies the relative signal strength for a pair of edges between $\maP$ and $\maQ$. 
Let $I$ denote the Kullback–Leibler (KL) divergence between a pair of correlated edges $(uv,\pi^*(u)\pi^*(v))$ under $\maP$ and $\maQ$:
\begin{align*}
    I = p_{11}\log \frac{p_{11}}{p^2}+2p_{10}\log \frac{p_{10}}{p(1-p)}+p_{00}\log\frac{p_{00}}{(1-p)^2}.
\end{align*}
Our goal is to determine the minimal sample size required to test the two hypotheses between $\maH_0$ and $\maH_1$.

It is well-known that the minimal value of the sum of Type I and Type II errors between ${\maP}$ and $\maQ$ is $1-\TV({\maP},{\maQ})$ (see, e.g., \cite[Theorem 13.1.1]{lehmann2005testing}), achieved by the likelihood ratio test, where $\TV({\maP},{\maQ}) = \frac{1}{2}\int |d{\maP}-d{\maQ}|$ is the total variation distance between ${\maP}$ and ${\maQ}$. Thus detection criterion~\eqref{eq:strong-detection} is equivalent to $\lim_{n\to\infty} \TV(\maP,\maQ) = 1$, while impossibility of testing is equivalent to \(\lim_{n\to\infty}\TV(\maP,\maQ)=0\). We first define the function
\[
C_{\varrho}(b)
:=
\begin{cases}
2, & b\leq 0,\\[2mm]
b\varrho^{-1}(b^{-1}), & b>0.
\end{cases}
\]
Here $\varrho$ is the densest subgraph function which will be introduced in~\eqref{eq:densest-subgraph}.
Let
\(
p=n^{-a+o(1)}\) and \( \gamma=n^{b+o(1)}
\) with $b\le a<1$,  and define $\mathcal R
\triangleq
\left\{b\leq0,\ a-b<\frac12\right\}
\cup
\left\{0<b\leq a,\ 2a-b<1\right\}$ (see Figure~\ref{fig:main} for details).
Throughout the paper, we focus on $(a,b)\in\mathcal R$. The complementary strict regimes are impossible even under full sampling (see Remark~\ref{rmk:full-sample} for details). Next, we introduce our main theorem.

\begin{theorem}\label{thm:main-IT}
Assume $p = n^{-a+o(1)}$ and $\gamma = n^{b+o(1)}$ with constants $(a,b)\in \mathcal R$.
Then, for every fixed constant
\(\delta>0\), when $n\to \infty$,
\[
\TV(\maP,\maQ)=
\begin{cases}
1-o(1),
&
\displaystyle
s^2\geq \bigl( C_{\varrho}(b)+\delta\bigr)
\frac{n\log n}
     {I},\\
o(1),
&
\displaystyle
s^2\leq \bigl( C_{\varrho}(b)-\delta\bigr)
\frac{n\log n}
     {I}.
\end{cases}
\]
\end{theorem}

In view of Theorem~\ref{thm:main-IT}, for the induced subgraphs $G_1,G_2$ sampled from $\mathbf{G}_1,\mathbf{G}_2$, the size of the common vertex set under the ground-truth permutation $\pi^*$ follows a hypergeometric distribution $\HG(n,s,s)$ (see Section~\ref{sec:info-density} for details), which concentrates around its expectation $s^2/n$. Here $\HG(N,K,M)$ denotes the law of
the number of successes in $M$ draws without replacement from a population of size $N$
containing $K$ successes.
The extra $\log n$ factor arises from a union bound over all possible latent mappings $\pi:V(G_1)\to V(G_2)$.
The factor $1/I$ originates from  the Kullback--Leibler (KL) divergence between the joint law of a correlated edge being present and the law of two independent edges being present. Combining these ingredients yields the sample-size requirement $s \gtrsim \sqrt{n\log n/I}\vee \sqrt{n} = \sqrt{n\log n/I}$.
For comparison, ~\cite{huang2025sample} studied the analogous problem in the Gaussian Wigner model. There, under $\mathcal H_0$ the edges are independent standard normals and under $\mathcal H_1$ corresponding edges are jointly normal with correlation $\rho$ under the latent permutation $\pi^*$; the optimal sample complexity is $s \asymp \sqrt{n\log n/\log(1/(1-\rho^2))}\vee \sqrt{n}$, which matches the above form with the KL-term $\log(1/(1-\rho^2))$ playing the role analogous to $I$ in the Bernoulli setting.

The constant $C_{\varrho}(b)$ is the sharp multiplicative factor in the
critical squared sample size in Theorem~\ref{thm:main-IT}. For $b\leq0$,
aggregate edgewise information determines the threshold and
$C_{\varrho}(b)=2$. For $b>0$, the densest subgraph of the intersection
graph becomes decisive and
$C_{\varrho}(b)=b\varrho^{-1}(b^{-1})$.
Proposition~\ref{prop:phase-transition} shows that
\[
    \lim_{b\downarrow0}b\varrho^{-1}(b^{-1})=2,
\]
so the transition at $b=0$ is continuous. Note that
\begin{align}\label{eq:asy-I}
    I = \begin{cases}
        n^{-2(a-b)+o(1)}, &\text{when }b\leq 0\\
        n^{-(2a-b)+o(1)}, &\text{when }b>0
    \end{cases},
\end{align} 
comparison with the constraint
$s\leq n$ yields the boundaries $a-b=1/2$ and $2a-b=1$, respectively;
the equality cases require a separate examination of the precise
lower-order terms.
Figure~\ref{fig:main} summarizes the resulting phase diagram. The
detectable interiors are given by $b\leq0$ and $a-b<1/2$, or by
$0<b\leq a$ and $2a-b<1$. Beyond these boundaries, detection is
impossible for every $s\leq n$. This full-sampling impossibility region
agrees with the thresholds for the fully observed model established
in~\cite{wu2023testing,ding2023detection}. Accordingly, our main focus
is on the two intermediate regimes in which induced-subgraph sampling
leads to a nontrivial sample-complexity threshold. Finally, the region
$b>a$ lies outside the model because it would require $\rho>1$.

Unlike the fully observed setting, our parameter range includes the nonvanishing-information regime $I=\Theta(1)$. Previous full-observation results focus on sparse regimes in which $I=o(1)$ and use  $H(p,\gamma)\triangleq p^2\big((1+\gamma)\log(1+\gamma)-\gamma\big)$
for the per-edge information~\cite{wu2023testing,ding2023detection}. In the present setting, the sharp threshold is governed by the exact mutual information $I$. When $p$ and $\rho$ converge to fixed nonzero constants, $I/H(p,\gamma)$ does not converge to one; hence replacing $I$ by $H(p,\gamma)$ would generally lose the sharp leading constant.

\begin{remark}[Full-sampling impossibility]\label{rmk:full-sample}
Let $\maP_n,\maQ_n$ denote the laws under full observation. Since
induced-subgraph sampling is a common Markov kernel under both
hypotheses,
\[
    \TV(\maP_s,\maQ_s)\le \TV(\maP_n,\maQ_n),\qquad s\le n.
\]
If $b\le0$ and $a-b>1/2$, then $\rho^2=o(n^{-1})$, so full-observation
impossibility follows from~\cite[Theorem~3]{wu2023testing}. If $b>0$
and $2a-b>1$, then $np_{11}=o(1)$, and the corresponding result follows
from~\cite[Theorem~1.1]{ding2023detection}. Hence detection is
impossible for every $s\le n$ in these two regions. 
Throughout the
paper, we therefore focus on the  regimes
$b\le0,\ a-b<1/2$ and $0<b\le a,\ 2a-b<1$.
\end{remark}

\begin{remark}
For $a>1$, the entire feasible range $b\le a$ lies in the region where detection is impossible even with two original graphs. The same holds for $a=1$ and $b<1$, leaving $(a,b)=(1,1)$ as the only boundary case. 
Indeed,  when $s^2 \ge (C_\varrho(b)+\delta) n\log n/I$ for $a = 1$, we still have  $\TV(\maP,\maQ) = 1-o(1)$. However, for the converse side, the total variation distance may remain bounded away from both zero and one, so the zero--one transition in Theorem~\ref{thm:main-IT} need not hold. This phenomenon has also been observed in the full-observation model; see \cite{ding2023detection,feng2025strong}. We therefore assume $b\le a<1$ throughout the paper.
\end{remark}



\begin{figure}[htbp]
    \centering
\begin{tikzpicture}[
    x=7.5cm,y=5cm,
    >=Latex,
    axis/.style={->,line width=0.55pt},
    bd/.style={line width=0.82pt,dashed},
    dom/.style={line width=0.9pt},
    lab/.style={font=\scriptsize,inner sep=1.4pt},
    smallbox/.style={font=\scriptsize,align=left,inner sep=2.6pt,rounded corners=1.3pt,fill=white,fill opacity=.90,text opacity=1},
    legendbox/.style={font=\scriptsize,align=left,inner sep=3.2pt,rounded corners=1.5pt,fill=white,fill opacity=.96,text opacity=1,draw=black!25}
  ]

\def\xmax{1.08}
\def\ymin{-0.56}
\def\ymax{1.06}

\fill[gray!13] (0,\ymin) -- (1,\ymin) -- (1,0) -- (0.5,0) -- (0,-0.5) -- cycle;
\fill[pattern=dots,pattern color=gray!75] (0,\ymin) -- (1,\ymin) -- (1,0) -- (0.5,0) -- (0,-0.5) -- cycle;
\fill[gray!13] (0.5,0) -- (1,0) -- (1,1) -- cycle;
\fill[pattern=dots,pattern color=gray!75] (0.5,0) -- (1,0) -- (1,1) -- cycle;

\fill[blue!10] (0,-0.5) -- (0,0) -- (0.5,0) -- cycle;
\fill[pattern=horizontal lines,pattern color=blue!63] (0,-0.5) -- (0,0) -- (0.5,0) -- cycle;

\fill[orange!14] (0,0) -- (1,1) -- (0.5,0) -- cycle;
\fill[pattern=north east lines,pattern color=orange!80!black] (0,0) -- (1,1) -- (0.5,0) -- cycle;

\fill[white,opacity=.73] (0,0) -- (0,\ymax) -- (1.02,\ymax) -- (1,1) -- cycle;

\draw[dom] (0,0) -- (1,1) node[pos=.60,above left=2pt,sloped,lab] {$b=a\;(\rho\le 1)$};
\draw[dom] (0,-0.5) -- (0.5,0) node[pos=.50,below=2pt,sloped,lab] {$a-b=\tfrac12$};
\draw[dom] (0.5,0) -- (1,1) node[pos=.58,below=3pt,sloped,lab] {$2a-b=1$};
\draw[dotted,black!45] (1,\ymin) -- (1,1.015);

\draw[axis] (0,\ymin) -- (0,\ymax) node[above,lab] {$b$};
\draw[axis] (0,0) -- (\xmax,0) node[right,lab] {$a$};
\foreach \x/\lbl in {0.5/{\tfrac12},1/{1}}{
  \draw (\x,0.011) -- (\x,-0.011) node[below=2pt,lab] {$\lbl$};
}
\foreach \y/\lbl in {-0.5/{-\tfrac12},0/{0},1/{1}}{
  \draw (0.011,\y) -- (-0.011,\y) node[left=2pt,lab] {$\lbl$};
}

\node[smallbox,anchor=center] at (0.145,-0.115)
  {$b\le0,\; a-b<\tfrac12$\\[-1pt]
   $s_*^2=\frac{2n\log n}{I}$};

\node[smallbox,anchor=center] at (0.36,0.14)
  {$0<b\le a,\; 2a-b<1$\\[-1pt]
   $s_*^2=\frac{b\varrho^{-1}(b^{-1})n\log n}{I}$};

\node[smallbox,anchor=center] at (0.72,-0.30)
  {full-sampling impossible\\[-1pt]
   $s\le n\Longrightarrow\TV(\maP,\maQ)=o(1)$};

\node[smallbox,anchor=center] at (0.36,0.83)
  {outside model\\[-1pt]
   $b>a$ would force $\rho>1$};

\node[legendbox,anchor=north west] at (1.05,0.97) {
\begin{tabular}{@{}c@{\ }l@{}}
\tikz{\fill[blue!10] (0,0) rectangle (.12,.08);\fill[pattern=horizontal lines,pattern color=blue!63] (0,0) rectangle (.12,.08);}
& constant $2$\\[-1pt]
\tikz{\fill[orange!14] (0,0) rectangle (.12,.08);\fill[pattern=north east lines,pattern color=orange!80!black] (0,0) rectangle (.12,.08);}
& constant $b\varrho^{-1}(b^{-1})$\\[-1pt]
\tikz{\fill[gray!13] (0,0) rectangle (.12,.08);\fill[pattern=dots,pattern color=gray!75] (0,0) rectangle (.12,.08);}
& impossible for all $s\le n$
\end{tabular}
};

\node[font=\scriptsize,align=center,anchor=south east] at (0.78,1.055)
  {$p=n^{-a+o(1)},\quad \gamma=n^{b+o(1)}$};

\end{tikzpicture}
    \caption{Phase diagram of the sharp sample-complexity threshold under $p=n^{-a+o(1)}$ and $\gamma=n^{b+o(1)}$. The blue and orange regions correspond to $C_{\varrho}(b)=2$ and $C_{\varrho}(b)=b\varrho^{-1}(b^{-1})$, respectively; in the gray region, detection is impossible for every $s\leq n$, while $b>a$ lies outside the model. The boundary cases are determined by checking the corresponding sample-size condition in Theorem~\ref{thm:main-IT} under the constraint $s\le n$.}
    \label{fig:main}
\end{figure}

\subsection{Related work}

\emph{Information-theoretic results.}
For the correlated \ER graphs model, the sharp threshold for dense graphs and the threshold within a constant factor for sparse graphs were determined in~\cite{wu2023testing}. It is shown in \cite{ding2023detection} that the constant for sparse graphs can be sharpened by analyzing the densest subgraphs. Another related model is the Gaussian Wigner model, where the weighted edges follow standard normals. The optimal detection threshold in the Gaussian Wigner was also derived in~\cite{wu2023testing}. More recently, \cite{huang2025sample} investigated the case where two induced subgraphs are sampled from the Gaussian Wigner model and derived the optimal sample complexity for correlation detection.

\emph{Polynomial-time algorithm and computational hardness.} There are several efficient algorithms for correlation detection in the correlated \ER graphs.
\cite{barak2019nearly} proposed a method based on counting balanced subgraphs, achieving weak detection criterion, where a balanced graph is the graph that is denser than any of its subgraphs. It is shown in~\cite{wu2023testing} that a polynomial-time algorithm based on counting trees succeeds for strong detection when $p\ge n^{-1+o(1)}$ and $\rho > \sqrt{\alpha}\approx \sqrt{0.338}$, where $\alpha$ is the Otter's constant introduced in~\cite{otter1948number}. The recent work~\cite{huang2025testing} considered counting bounded degree motifs and achieved weak detection under the conditions $p\ge n^{-2/3}$ and  $\rho = \Omega(1)$.

Inspired by the sum-of-squares method, the low-degree framework provides a powerful method for deriving computational lower bounds in a wide range of high-dimensional statistical problems~\citep{hopkins2017power,hopkins2018statistical}. 
\cite{mao2024testing} formulated a low-degree framework suggesting that polynomial-time detection should not be possible when $\rho^2 \le 1/\mathrm{polylog}(n)$. 
Further evidence is provided in~\cite{ding2023low}, which shows that  any degree-$O(\rho^{-1})$ polynomial-time algorithm fails for strong detection. Moreover, it is conjectured that any degree-$D$ polynomial-time algorithm fails for detection as long as $\log D  = o\pth{\frac{\log n}{\log np}\wedge \sqrt{\log n}}$ and $\rho<\sqrt{\alpha}$ in the sparse regime where $p = n^{-1+o(1)}$.

\emph{Graph matching.}
The graph matching problem refers to finding a correspondence between the nodes of different graphs~\citep{caetano2007learning}. Recently, substantial progress has been made on matching two correlated random graphs, especially the correlated \ER model and the correlated Gaussian Wigner model. The statistical analysis includes~\cite{cullina2016improved,cullina2017exact,cullina2020partial,ganassali2021impossibility,wu2022settling,ding2023matching,du2025optimal}. Among these studies, \cite{wu2022settling} determined the information-theoretic thresholds for exact recovery and characterized the thresholds for partial recovery up to a constant. The constant gap was filled by~\cite{ding2023matching} by analyzing the densest graph.
Beyond information-theoretic analysis, there are also many efficient algorithms proposed for graph matching, including~\cite{shirani2017seeded,barak2019nearly,dai2019analysis,ganassali2020tree,mossel2020seeded,ding2021efficient,ding2023polynomial,mao2023exact,fan2019spectral,ganassali2024statistical,ganassali2024correlation,ding2025polynomial,mao2025random}. There are also many extensions on Gaussian Wigner model and correlated \ER graph model, including the inhomogeneous \ER model~\citep{racz2023matching, ding2023efficiently}, the partially correlated graphs model~\citep{huang2024information}, the correlated stochastic block model~\citep{chen2024computational,chen2025detecting}, the multiple correlated graphs model~\citep{ameen2024robust,ameen2025detecting}, the correlated random geometric graphs model~\citep{wang2022random}, and the attributed \ER model~\citep{zhang2024attributed,yang2024exact}.
\subsection{Notation}
For $n\in\mathbb{N}$,
let $[n]\triangleq\{1,\ldots,n\}$. We use the standard asymptotic
notation $O(\cdot)$, $\Omega(\cdot)$, $o(\cdot)$, $\omega(\cdot)$,
$\lesssim$, $\gtrsim$, and $\asymp$, where the implicit constants are
independent of $n$.
For a graph $G$, let $V(G)$ and $E(G)$ denote its vertex and edge sets,
and write $|V(G)|$ and $|E(G)|$ as the sizes. We write
$uv$ for the unordered pair $\{u,v\}$ and let
\(
    \binom{S}{2}\triangleq\{uv:u,v\in S,\ u\neq v\}.
\)
For an injective mapping $\pi:S\to T$ and an edge
$e=uv\in\binom{S}{2}$, define $\pi(e)\triangleq\pi(u)\pi(v)$. In our sampling model, $\mathbf{G}_1,\mathbf{G}_2$ denote the original
$n$-vertex graphs, whereas $G_1,G_2$ denote the sampled $s$-vertex
induced subgraphs. We use $\mathcal{Q}$ and $\mathcal{P}$ for their
joint laws under $\mathcal{H}_0$ and $\mathcal{H}_1$, respectively.

\section{Information density between two correlated graphs}\label{sec:info-density}

In this section, we introduce some basic properties between the correlated subgraphs $G_1$ and $G_2$ under our sampling model. For the induced subgraphs $G_1,G_2$ sampled from $\mathbf{G}_1,\mathbf{G}_2$, we denote the set of common vertices as 
\begin{align}\label{eq:def_of_Spi*}
    \sfS_{\pi} = V(G_1)\cap \pi^{-1}(V(G_2)),\quad \sfT_{\pi} = \pi(V(G_1))\cap V(G_2).
\end{align}
Under both $\maH_0$ and $\maH_1$, $|\sfS_\pi|=|\sfT_\pi|$ follow the Hypergeometric distribution $\HG(n,s,s)$ with mean $s^2/n$ and variance $\frac{s^2}{n}\cdot \frac{(n-s)^2}{n(n-1)}\le \frac{s^2}{n}$. By Chebyshev's inequality,
\[
\prob{\left||\sfS_\pi|- \frac{s^2}{n}\right|\ge \epsilon \frac{s^2}{n}}
\le \frac{\var{|\sfS_\pi|}}{\epsilon^2 (s^2/n)^2}
\le \frac{1}{\epsilon^2 (s^2/n)}.
\]
Thus, for any fixed $\epsilon>0$, if $s^2/n=\omega(1)$, then $|\sfS_\pi|=|\sfT_\pi|\in R_\epsilon\triangleq [(1-\epsilon)s^2/n,(1+\epsilon)s^2/n]$ with probability $1-o(1)$. For any $r\in\mathbb N$, denote by $\maP_{r}$ the distribution conditional on $|\sfS_\pi|=|\sfT_\pi|=r$, then
\begin{align}\label{eq:epsilon}
    \TV(\maP,\maQ)\le o(1)+\sum_{r\in R_\epsilon}\prob{|\sfS_\pi|=r}\TV(\maP_r,\maQ), \text{ where } R_\epsilon=\Big[\frac{(1-\epsilon)s^2}{n},\frac{(1+\epsilon)s^2}{n}\Big]\cap\mathbb N.
\end{align}
Consequently, it suffices to prove $\sup_{r\in R_\epsilon}\TV(\maP_{r},\maQ)=o(1)$ for impossibility results. The constant $\epsilon$ in~\eqref{eq:epsilon} will be chosen later.
Under $\maP_{r}$, there is a latent injective mapping $\pi:\sfS\subseteq V(G_1)\to \sfT\subseteq V(G_2)$, where $\pi$ is uniformly distributed over all such mappings with $|\sfS|=|\sfT|=r$. We denote $\Pi_r$ the set of such injective mappings.
We assume $r\in R_\epsilon$ throughout the paper.

We then introduce the information density under $\maP_r$. Let $G(e)\triangleq \indc{e\in E(G)}$. Given two graphs $G_1$ and $G_2$, we define an \emph{information graph} with vertex set $\sfS_\pi$ and weighted edge $i_\pi(e)\triangleq \log\tfrac{\maP_r(G_1(e),G_2(\pi(e)))}{\maP_r(G_1(e))\maP_r(G_2(\pi(e)))}$, where $\pi:\sfS_\pi\mapsto \sfT_\pi$ is latent permutation (see also the information density defined in~\cite[Section 18.1]{polyanskiy2025information}). Denote $\maI_\pi(V)$ as the set of edges $i_\pi(uv)$ for $u,v\in V\subseteq \sfS_\pi$ and $|\maI_\pi(V)| = \sum_{\{u,v\}\subseteq V}i_\pi(uv)$ as the total weights. Recall that $I = \E_{\maP_r}[i_\pi(e)]$. 
By~\cite{anantharam2016densest}, for any \ER graph $G\sim \maG(n,\tfrac{\lambda}{n})$ with  any constant $\lambda>0$, there exists a constant $\varrho(\lambda)>0$ such that 
\begin{align}\label{eq:densest-subgraph}
    \max_{\emptyset\neq U\subseteq V(G)} \frac{|\maE(U)|}{|U|}\overset{\P}{\longrightarrow } \varrho(\lambda )\text{ as }n\to\infty, 
\end{align}
where $\maE(U)$ denotes the collection of edge in vertex sets $U$. Moreover, $\varrho$ is increasing and continuous, and it follows from~\cite[Proposition 2.3]{ding2023detection} that one can define the inverse function $\varrho^{-1}:[1,\infty)\to [1,\infty)$ with $\varrho^{-1}(1) = 1$.
We now establish the following proposition on information density.

\begin{proposition}\label{prop:information-density}
If $rI = \omega(1\vee \log(1+\gamma))$, then \begin{align}\label{eq:info-density-bounded}
    \max_{V\neq \emptyset} \frac{|\maI_\pi(V)|}{|V|} = \frac{rI}{2}(1+o_p(1)),\qquad \text{as } r\to\infty.
\end{align}    
If $\gamma = \omega(1)$ and $\frac{rI}{\log(1+\gamma)}\to c$ for some constant $c>0$, then \begin{align}\label{eq:info-density-poly}
    \max_{V\neq \emptyset} \frac{|\maE_\pi(V)|}{|V|} = \varrho(c) (1+o_p(1)),\qquad \text{as } r\to\infty,
\end{align}
where $\varrho:(0,\infty)\mapsto (0,\infty)$ is a continuous and increasing function.
\end{proposition}

\begin{proof}
We first focus on the case $r I = \omega(1\vee \log(1+\gamma))$.
Note that $\maQ(G_1(e), G_2(\pi(e)))=\maP_r(G_1(e)) \maP_r( G_2(\pi(e)) )$.
Define $R\triangleq \frac{\diff \maP_r}{\diff \maQ}(G_1(e), G_2(\pi(e)))$.
Then $R\le 1+\gamma$. Let $B\triangleq \log(1+\gamma)$.
For one edge, it follows from Lemma~\ref{lem:var-gamma-bounded} that $\var_{\maP_r}[i_\pi(e)]\le (2+B)I.$
For a fixed subset $V$ with $|V|=k$, by Chebyshev's inequality, for any constant $\delta>0$,
\[
\prob{|\calI_\pi(V)| \le \Expect|\calI_\pi(V)| - \delta r k I }
\le \frac{\frac{k^2}{2}(2+B) I}{(\delta r k I)^2}\to 0.
\]
Therefore, $\frac{|\calI_\pi(V)|}{|V|}
\ge \frac{rI}{2}(1+o_p(1))$ for $|V|=r$.
Additionally, by Bennett's inequality (see, e.g., ~\cite[Theorem 2.9]{boucheron2013concentration}),
\[
\prob{|\calI_\pi(V)| \ge \Expect|\calI_\pi(V)| + \delta r k I }
\le \exp\pth{- k \frac{ \delta  r I }{B }\phi\pth{\frac{B \delta r k I}{\frac{k^2}{2}(2+B)I }}},
\]
where $\phi(x)=\frac{(1+x)\log(1+x)-x}{x}$.
Using $\phi(x)\gtrsim x$ for $x\leq 1$ and $\phi(x)\gtrsim 1+\log x$ for $x>1$, we have uniformly over $1\leq k\leq r$ that
\(
\frac{rI}{B}\phi\left(\frac{2B\delta}{2+B}\frac{r}{k}\right)
\gtrsim
\frac{rI}{1+B}\log\left(\frac{er}{k}\right).
\)
Thus, the condition $rI \gg 1\vee B$ implies that $\frac{rI}{B}\phi(\frac{2B\delta}{2+B}\frac{r}{k})\gg \log(\frac{er}{k})$.
Applying the union bound over $|V|=k$ with $\binom{r}{k}\le (\frac{er}{k})^k$ and then over $k\in[r]$, we obtain
\[
\prob{\exists V\ne \emptyset: |\calI_\pi(V)| \ge \Expect|\calI_\pi(V)| + \delta r k I }
\le \sum_{k=1}^r \exp\pth{-2k\log(er/k)} 
\to 0.
\]
The upper bound $\max_{V\ne \emptyset} \frac{|\calI_\pi(V)|}{|V|}
\le \frac{rI}{2}(1+o_p(1))$ follows since $\Expect|\calI_\pi(V)|\le \frac{rk}{2}I$.

We then focus on the case $\gamma = \omega(1)$ and $\tfrac{rI}{\log(1+\gamma)} \to c$ for some constant $c>0$. We have $I = (1+o(1))p_{11}\log(1+\gamma)$ when $\gamma = \omega(1)$, and thus $r p_{11}\to c$. Given $G_1$ and $G_2$, define the intersection graph with vertex set $\sfS_\pi$ and edge $G_1(e)\cap G_2(\pi(e))$. Denote $\maE_\pi(V)$ as the set of edges $uv$ for $u,v\in V\subseteq \sfS_\pi$ and $|\maE_\pi(V)|$ as the total edges. It follows from~\cite{anantharam2016densest} that the densest subgraph satisfies \begin{align}\label{eq:densest-func}
    \max_{\emptyset \neq V \subseteq \sfS_\pi} \frac{|\maE_\pi(V)|}{|V|} = \varrho(c) (1+o_p(1)),\qquad\text{ as }r\to\infty,
\end{align}
where $\varrho(x)$ is a continuous and increasing function.
\end{proof}

\section{Information-theoretic lower bounds}

In this section, we introduce the key points on information-theoretic lower bounds, under which the total variation distance $\TV(\maP_r,\maQ) = o(1)$. It is well known that the total variation distance is upper bounded by the second moment (see, e.g., \cite[Equation 2.27]{tsybakov2009introduction}): $\TV(\maP_r,\maQ)\le \sqrt{\chi^2(\maP_{r}\Vert \maQ)}$,
where $\chi^2(\maP_{r}\Vert \maQ)\triangleq \E_{\maQ}\big[\big(\tfrac{\maP_{r}}{\maQ}-1\big)^2\big]$.
However, the second moment sometimes explodes due to certain rare events, while $\TV({\maQ},{\maP_{r}})$ remains bounded away from one. In order to circumvent such catastrophic events, one method is using the conditional second moment method on these events. Specifically, given some high probability event $\maE$ with ${\maP_{r}}(\maE) = 1-o(1)$, we define the conditional distribution as $\maP_{r}' = \maP_{r}'(\cdot|\maE)$.
A sufficient condition for establishing impossibility via the conditional second moment method is
\begin{align}\label{eq:condi_second_moment}
    \chi^2(\maP_{r}'\Vert \maQ) = o(1)
    \quad\Longrightarrow\quad
    \TV(\maP_r,\maQ) = o(1).
\end{align}

\subsection{Necessity of conditional second moment}

We have shown in~\eqref{eq:condi_second_moment} that the conditional second moment method suffices to control the total variation distance, provided that one can construct a suitable event $\maE$ with $\maP_r(\maE) = 1-o(1)$. However, finding such an event may be challenging. This raises a natural question: is this conditioning step in some sense necessary, or at least equivalent to the assumption that $\TV(\maP,\maQ)=o(1)$? Equivalently, one may ask:

\emph{Does $\TV(\maP,\maQ)=o(1)$ imply the existence of an event $\maE$ such that $\chi^2(\maP(\cdot\mid \maE)\|\maQ)=o(1)$?}

In this subsection, we provide a construction of such event which provides a positive answer to this question. 
\begin{lemma}\label{lem:general-condi-second-moment}
    Let $\maP$ and $\maQ$ be probability measures on the same measurable space, and let $\tau = \TV(\maP,\maQ)\le 0.1$, then \begin{align*}
        \inf_{\maE}\pth{\chi^2(\maP(\cdot\mid\maE)\Vert \maQ)+\maP(\maE^c)}\le 5\tau.
    \end{align*}
\end{lemma}
\begin{proof}
    Let $\maP = \maP_{ac}+\maP_\perp$ be the Lebesgue decomposition of $\maP$ with respect to $\maQ$. Here $\maP_{ac} \ll\mathcal Q$ is the absolutely continuous
part of $\mathcal P$, whereas
$\mathcal P_{\perp}\perp\mathcal Q$ is its singular part. Define the measurable space as $\Omega$. Choose a measurable set $A$ such that $\maQ(A) =0$ and $\maP_\perp(A^c) = 0$.
Let $L\triangleq \frac{d\,\maP_{ac}}{d\,\maQ}$.
Since $\TV(\maP,\maQ) = \tau$, we have \begin{align*}
    \E_\maQ[(L-1)_+]+\maP_\perp(\Omega) = \tau.
\end{align*} 
Let $\maE_0 = \sth{L\le 3}\cap A^c$. Since \(L\le\frac32(L-1)\) on \(\{L>3\}\),
\begin{align*}
\mathcal P(\mathcal E_0^c)
&=\mathcal P_{\perp}(\Omega)
  +\mathbb E_{\mathcal Q}
     \left[L\mathbf 1_{\{L>3\}}\right]\le
\mathcal P_{\perp}(\Omega)
+\frac32\mathbb E_{\mathcal Q}
  \left[(L-1)\mathbf 1_{\{L>3\}}\right]\le
\frac32\left(
\tau-
\mathbb E_{\mathcal Q}
\left[(L-1)\mathbf 1_{\{1<L\le3\}}\right]
\right).
\end{align*}
In particular, $\mathcal P(\mathcal E_0)
\ge1-\frac32\tau>0.$ Since $\maE_0$ excludes the singular part, we have $\frac{d\,\maP(\cdot \mid \maE_0)}{d\,\maQ} = \frac{L\indc{\maE_0}}{\maP(\maE_0)}$. Since $L^2\le L$ when $L\le 1$ and $L^2\le L+3(L-1)$ when $1\le L\le 3$, we have \begin{align*}
    \mathbb E_{\mathcal Q}
\left[L^2\mathbf 1_{\mathcal E_0}\right]
\le
\mathcal P(\mathcal E_0)
+
3\mathbb E_{\mathcal Q}
\left[(L-1)\mathbf 1_{\{1<L\le3\}}\right].
\end{align*}
It follows that \begin{align*}
    &~\chi^2(\maP(\cdot \mid\maE_0)\Vert \maQ)+\maP(\maE_0^c) = \frac{\E_\maQ[L^2\indc{\maE_0}]}{\maP(\maE_0)^2}-1+\maP(\maE_0^c)\\\le&~ \frac{\maP(\maE_0)-\maP(\maE_0)^3+3\E_\maQ[(L-1)\indc{1\le L\le 3}]}{\maP(\maE_0)^2}\overset{\mathrm{(a)}}{\le} \frac{2\maP(\maE_0^c)+3\E_\maQ[(L-1)\indc{1\le L\le 3}]}{\maP(\maE_0)^2}\overset{\mathrm{(b)}}{\le} 5\tau,
\end{align*}
where $\mathrm{(a)}$ is because $x-x^3 = x(1-x)(1+x)\le 2(1-x)$ with $x = \maP(\maE_0)\in [0,1]$; $\mathrm{(b)}$ follows from  $2\maP(\maE_0^c)+3\E_\maQ[(L-1)\indc{1\le L\le 3}]\le 3\TV(\maP,\maQ) = 3\tau$, $\maP(\maE_0)\ge 1-\tfrac{3\tau}{2}$, and $\tfrac{3\tau}{(1-3\tau/2)^2}\le 5\tau$ when $\tau\le 0.1.$
\end{proof}

By Lemma~\ref{lem:general-condi-second-moment}, the conditional second moment method incurs no loss of information. We will therefore employ this method in the subsequent information-theoretic lower bound analysis. Furthermore, the lemma suggests that truncation based on the likelihood ratio is a natural and principled choice.

\subsection{Core set and lower bound analysis for weak signal}\label{subsec:weak-signal-lwbd}

We next prove the sharp lower bound in Theorem~\ref{thm:main-IT} for $b\le0$. Recall that $\Pi_r$ is the collection of injective mappings from $\sfS_\pi\subseteq V(G_1)$ to $\sfT_\pi\subseteq V(G_2)$ with $|\sfS_\pi|=|\sfT_\pi|=r$. We write $\pi\bot\ti\pi$ when $\pi$ and $\ti\pi$ are sampled independently and uniformly from $\Pi_r$, and put
\[
L_\pi\triangleq\prod_{e\in\binom{\sfS_\pi}{2}}
\ell\bigl(G_1(e),G_2(\pi(e))),
\]
where $\ell:\{0,1\}^2\to\mathbb{R}_{+}$ is the single-edge likelihood ratio defined by
\[
\ell(x,y)\triangleq
\frac{\maP_r\bigl(G_1(e)=x,G_2(\pi(e))=y\mid\pi\bigr)}
{\maQ\bigl(G_1(e)=x,G_2(\pi(e))=y\bigr)},
\qquad x,y\in\{0,1\}.
\]
For an event $\maE$ with $\maP_r(\maE)=1-o(1)$, define $\maP_r'=\maP_r(\cdot\mid\maE)$ and
\begin{align}\label{eq:def-psi-weak}
\psi(\pi,\ti\pi)
\triangleq
\log\E_{\maQ}\left[L_\pi L_{\ti\pi}
\indc{\maE_\pi}\indc{\maE_{\ti\pi}}\right].
\end{align}
Then
\begin{align}\label{eq:conditional-second-moment-psi}
1+\chi^2(\maP_r'\Vert\maQ)
=\frac{\E_{\pi\bot\ti\pi}[e^{\psi(\pi,\ti\pi)}]}{\maP_r(\maE)^2}.
\end{align}
It therefore suffices to prove $\E_{\pi\bot\ti\pi}[e^{\psi(\pi,\ti\pi)}]=1+o(1)$.

We first consider a structural decomposition of $L_\pi L_{\ti\pi}$. The independent part can be analyzed via the correlated functional digraph associated with $(\pi,\ti\pi)$ (see Remark~\ref{rmk:correlated-func-digraph} for details). Since this graph has maximum degree two, each of its connected components is either a path or a cycle.
Let $\sfP^\sfE$ and $\sfC^\sfE$ denote the collections of edge paths and closed edge cycles, respectively. If $P=(e_1,\pi(e_1),\ldots,e_j,\pi(e_j))$ is a path, then $\ti\pi(e_{i+1})=\pi(e_i)$ for $1\le i<j$. A cycle satisfies the additional relation $\ti\pi(e_1)=\pi(e_j)$. For $P\in\sfP^\sfE$ and $C\in\sfC^\sfE$, let
\begin{align*}
L_P&\triangleq
\prod_{e\in\binom{\sfS_\pi}{2}\cap P}\ell\bigl(G_1(e),G_2(\pi(e))\bigr)
\prod_{e\in\binom{\sfS_{\ti\pi}}{2}\cap P}\ell\bigl(G_1(e),G_2(\ti\pi(e))\bigr),\\
L_C&\triangleq
\prod_{e\in\binom{\sfS_\pi}{2}\cap C}\ell\bigl(G_1(e),G_2(\pi(e))\bigr)
\prod_{e\in\binom{\sfS_{\ti\pi}}{2}\cap C}\ell\bigl(G_1(e),G_2(\ti\pi(e))\bigr).
\end{align*}
Thus $L_\pi L_{\ti\pi}=\prod_{P\in\sfP^\sfE}L_P\prod_{C\in\sfC^\sfE}L_C$. The following orbit calculation is standard; a proof is given in Appendix~\ref{apd:proof-lem-orbit-expect}.
\begin{lemma}\label{lem:orbit-expect}
For every $P\in\sfP^\sfE$ and $C\in\sfC^\sfE$,
\[
\E_{\maQ}[L_P]=1,
\qquad
\E_{\maQ}[L_C]=1+\rho^{2|C|}.
\]
\end{lemma}

The closed edge orbits are generated by a simpler collection of vertex cycles. This motivates the following definition.
\begin{definition}[Core set]\label{def:core-set}
The \emph{core set} of $(\pi,\ti\pi)$ is the unique largest invariant subset
\begin{align}\label{eq:def_of_I}
S^* = S^*(\pi,\ti\pi)
\triangleq
\bigcup\left\{S\subseteq\sfS_\pi\cap\sfS_{\ti\pi}:\pi(S)=\ti\pi(S)\right\}.
\end{align}
Equivalently, $S^*$ is the union of the vertex cycles of $\ti\pi^{-1}\circ\pi$.
\end{definition}
This definition includes isolated fixed points, which are not incident to an edge cycle when $|S^*|=1$. 
The closed edge cycles correspond to the orbits of the induced action of \(\tilde\pi^{-1}\circ\pi\) on \(\binom{S^*}{2}\), and these orbits partition \(\binom{S^*}{2}\).
In particular, restricting to $S^*$ reduces the cyclic part of the second moment to the fully correlated setting. Let $C_k$ be the number of vertex cycles of length $k$. Then
\begin{align}\label{eq:core-vertex-cycle-decomposition}
|S^*|=\sum_{k\ge1}kC_k.
\end{align}
The factorial moments of $(C_k)_{k\ge1}$ satisfy the following bound.
\begin{lemma}\label{lem:property_of_I}
For every sequence $(a_k)_{k\ge1}$ of nonnegative integers,
\begin{align}\label{eq:cycle-moment}
\E_{\pi\bot\ti\pi}\left[\prod_{k\ge1}(C_k)_{a_k}\right]
\le
\prod_{k\ge1}\Big(\frac{(r^2/s^2)^k}{k}\Big)^{a_k},
\end{align}
where $(x)_a=x(x-1)\cdots(x-a+1)$. Consequently,
\begin{align}\label{eq:S^*-bound}
\prob{C_1=c_1}\le\frac{(r^2/s^2)^{c_1}}{c_1!},
\qquad
\prob{|S^*|=t}\le(r^2/s^2)^t.
\end{align}
\end{lemma}
The proof is deferred to Appendix~\ref{apd:proof-lem-property-I}. This is a sub-Poisson-type bound in the factorial-moment sense: the
joint factorial moments of $(C_k)_{k\geq1}$ are dominated by those of
independent Poisson random variables with means
$\lambda_k=(r^2/s^2)^k/k$.
We then introduce our event. Define
\begin{align*}
    \maE_\pi^1\triangleq \sth{(G_1,G_2,\pi):\max_{V\neq \emptyset} \frac{|\maI_\pi(V)|}{|V|}\le \frac{rI}{2}\Big(1+\frac{\delta}{10}\Big)}
\end{align*}
It follows from Proposition~\ref{prop:information-density} that $\maP_r(\maE_\pi^1) = 1-o(1)$ if $rI = \omega(1\vee \log(1+\gamma))$.  The unconditional orbit bound and the truncation bound on $\maE_\pi^1$ give, respectively,
\[
\prod_{C\in\sfC^\sfE}(1+\rho^{2|C|})
\le(1+\rho^2)^{\binom{|S^*|}{2}}
\quad\text{and}\quad
\prod_{e\in\binom{S^*}{2}}\ell\bigl(G_1(e),G_2(\pi(e))\bigr)
\le \exp\Big(\frac{rI}{2}(1+\frac{\delta}{10})|S^*|\Big).
\]
By Lemma~\ref{lem:orbit-expect}, the remaining likelihood on path have expectation one under $\maQ$.
Let $\vp(t)\triangleq \min(\tfrac{rI}{2}(1+\tfrac{\delta}{10})t,\binom{t}{2}\log(1+\rho^2))$.
Hence
\begin{align}\label{eq:psi-core-envelope}
\exp\pth{\psi(\pi,\ti\pi)}
\le
\exp\Big(\varphi(|S^*|)\Big).
\end{align}
Moreover, when $s^2\le \tfrac{(2-\delta)n\log n}{I}$, we have $\tfrac{rI}{2}(1+\tfrac{\delta}{10})\le (1+\epsilon)(1-\tfrac{\delta}{2})(1+\tfrac{\delta}{10})\log n$. Pick $\epsilon$ sufficient small yields \begin{align}\label{eq:-delta4}
    \vp(t)\le \min\Big((1-\frac{\delta}{4})t\log n,\binom{t}{2}\log(1+\rho^2)\Big)\quad\text{ when }\quad s^2\le \frac{(2-\delta)n\log n}{I}.
\end{align}
We first present the bound obtained by using only the size of the core set.
\begin{proposition}\label{prop:core-size-mgf}
We have  
    \[
\E[\exp(\varphi(|S^*|))]
\le 1+ \sum_{k=2}^r \exp\pth{\varphi(k) + 2k \log \frac{r}{s} }.
\]
Moreover, if $s^2\le \tfrac{(2-\delta)n\log n}{I}$ and $|\log I| = o(\log n)$, then $\E[\exp(\varphi(|S^*|))] = 1+o(1)$.
\end{proposition}
\begin{proof}
    By Lemma~\ref{lem:property_of_I},\begin{align*}
        \E[\exp(\varphi(|S^*|))]&\le 1+\sum_{k=2}^r \exp(\vp(k))\prob{|S^*| = k}\le 1+\sum_{k=2}^r\exp (\vp(k) )(r^2/s^2)^k.
    \end{align*} 
    Since $\vp(t)\le (1-\tfrac{\delta}{4})t\log n$ and $\log(r^2/s^2)\le \log(2\log n/nI) = (-1+o(1))\log n$, we have \begin{align*}
        \E[\exp(\varphi(|S^*|))]\le 1+\sum_{k=2}^r\exp\Big((1-\frac{\delta}{4})k\log n+(-1+o(1))k\log n\Big) = 1+o(1).
    \end{align*}
\end{proof}

Proposition~\ref{prop:core-size-mgf} suffices for the desired lower bound under the condition $|\log I|=o(\log n)$. 
However, if $I$ is polynomially small, however, the
estimate $\prob{|S^*|=k}\le (r^2/s^2)^k$ treats all vertex cycles in the same way
and loses the sharp constant. We then isolate the fixed points in the following.

\begin{proposition}\label{prop:AFWnd-mgf}
We have 
    \[
\E[\exp(\varphi(C_1))]
\le 1+ \sum_{k=2}^r (\vp(k)\wedge 1)\exp\pth{\varphi(k) + k \log \frac{e r^2}{k s^2} }.
\]
Moreover, if $s^2\le \tfrac{(2-\delta)n\log n}{I}$ and $r\rho^2 \le n^{o(1)}$, then $\E[\exp(\vp(C_1))] = 1+o(1)$. 
\end{proposition}

\begin{proof}
Since $\varphi(0)=\varphi(1)=0$, the fixed-point bound
$\prob{C_1=k}\le (r^2/s^2)^k/k!$ gives
\begin{align}
\nonumber &~\E[\exp(\varphi(C_1))]=
1+\sum_{k=2}^r\bigl(\exp(\varphi(k))-1\bigr)
\prob{C_1=k}\\\label{eq:prop2-eq1}
\overset{\mathrm{(a)}}{\le}&~ 1+\sum_{k=2}^r(\exp(\varphi(k))-1)
\frac{(r^2/s^2)^k}{k!}\overset{\mathrm{(b)}}{\le} 1+\sum_{k=2}^r(\vp(k)\wedge 1)
\exp\pth{\varphi(k)+k\log\frac{er^2}{ks^2}},
\end{align}
where $\mathrm{(a)}$ follows from~\eqref{eq:S^*-bound}; $\mathrm{(b)}$ is because $k!\ge(k/e)^k$ and $e^x-1\le (x\wedge 1)e^x$ for $x\ge 0$.
We next prove the second claim. 
The claim is immediate when $\rho=0$, so assume $\rho>0$.

For $2\le k\le 1/\log(1+\rho^2)$, we have
\[
\varphi(k)
\le \binom{k}{2}\log(1+\rho^2)
\le \frac{k}{2}.
\]
Since $r^2/s^2\le1$ and
$r^2/s^2\le(1+\epsilon)r/n$, we have 
\begin{align}
\nonumber &~\sum_{k=2}^{1/\log(1+\rho^2)}
\vp(k)
\exp\pth{\varphi(k)+k\log\frac{er^2}{ks^2}}\le
\sum_{k\ge2}\frac{k^2\log(1+\rho^2)}{2}\exp\pth{\frac{k}{2}+k\log\frac{er^2}{ks^2}}\\\label{eq:prop2-eq2}\le&~ \frac{1}{2}\log(1+\rho^2)\sum_{k\ge 2}k^2\pth{\frac{\sqrt{e}r^2}{2s^2}}^k\le C\log (1+\rho^2)\left(\frac{r^2}{s^2}\right)^2
\le \frac{C(1+\epsilon)r\rho^2}{n}=o(1).
\end{align}
For $k>1/\log(1+\rho^2)$, we have 
\begin{align}
    \nonumber&~\sum_{k>1/\log(1+\rho^2)} \exp\pth{\vp(k)+k\log\frac{er^2}{ks^2}}\overset{\mathrm{(a)}}{\le}\sum_{k>1/\log(1+\rho^2)} n^{(1-\delta/4)k} \big(\frac{er^2}{ks^2}\big)^k\\\label{eq:prop2-eq3}\overset{\mathrm{(b)}}{\le}&~\sum_{k>1/\log(1+\rho^2)}\left(e(1+\epsilon)n^{-\delta/4}r\log(1+\rho^2)\right)^k\overset{\mathrm{(c)}}{\le} \sum_{k>1/\log(1+\rho^2)} n^{(-\delta/4-o(1))k} = o(1),
\end{align}
where $\mathrm{(a)}$ follows from~\eqref{eq:-delta4}; $\mathrm{(b)}$ is because $k>\tfrac{1}{\log(1+\rho^2)}$ and $\tfrac{r^2}{s^2}\le \tfrac{(1+\epsilon)r}{n}$; $\mathrm{(c)}$ follows from $r\log(1+\rho^2)\le r\rho^2\le n^{o(1)}$.
By~\eqref{eq:prop2-eq1},~\eqref{eq:prop2-eq2}, and~\eqref{eq:prop2-eq3}, we obtain $\E[\exp(\varphi(C_1))]=1+o(1)$.
\end{proof}

Proposition~\ref{prop:AFWnd-mgf} controls the closed edge cycles generated by pairs of fixed
points. It remains to show that all other closed edge cycles are negligible.
The following proposition provides this reduction.

\begin{proposition}\label{prop:AFWnd}
If $s^2\le \tfrac{(2-\delta)n\log n}{I}$ and  $r\rho^4=o(1)$, then
\[
\E[\exp(\psi(\pi,\ti\pi))]
\le \E[\exp(\varphi(C_1))]+o(1).
\]
\end{proposition}

The proof is deferred to Appendix~\ref{apd:proof-prop-AFWnd}.
It separates the remaining closed edge cycles into three classes: those
within one  vertex cycle, those joining a fixed point to a
 vertex cycle, and those joining two distinct  vertex
cycles. The condition $r\rho^4=o(1)$ makes every class negligible.
We are now ready to complete the lower-bound argument in the
subpolynomial-signal regime.

\begin{proof}[Proof of Theorem~\ref{thm:main-IT}, lower bound when
$b\le 0$]
Fix $\delta>0$ and suppose that
\(
s^2\le\tfrac{(2-\delta)n\log n}{I}.
\)
Let $\mathcal P(s)$ and $\mathcal Q(s)$ denote the joint laws of the two observed $s$-vertex induced subgraphs under the alternative and null hypotheses, respectively. For $s_1\le s_2$, independently retaining $s_1$ uniformly chosen vertices from each $s_2$-vertex observation yields the experiment with sample size $s_1$ under either hypothesis. Hence, by the data-processing inequality,
\begin{align}\label{eq:monoto-data-processing}
    \TV\bigl(\mathcal P(s_1),\mathcal Q(s_1)\bigr)
\le
\TV\bigl(\mathcal P(s_2),\mathcal Q(s_2)\bigr).
\end{align}
Therefore, it suffices to
consider the boundary case
\(
s^2=\tfrac{(2-\delta)n\log n}{I}.
\)
Hence, uniformly over
$r\in R_\epsilon$,
\begin{align}\label{eq:info-condition-dense}
    rI=\Theta(\log n)
=\omega\bigl(1\vee\log(1+\gamma)\bigr)
\end{align}
because $\gamma\le n^{o(1)}$. 
Then
Proposition~\ref{prop:information-density} gives
$\maP_r(\maE_\pi^1)=1-o(1)$
uniformly over $r\in R_\epsilon$.

By~\eqref{eq:asy-I}, we have $I = n^{-2(a-b)+o(1)}$ when $b\le 0$.
 We first consider $a=b$, or
equivalently $|\log I|=o(\log n)$. 
Proposition~\ref{prop:core-size-mgf} and the core-set bound therefore give
\[
\E_{\pi\bot\ti\pi}[\exp(\psi(\pi,\ti\pi))]
\le\E_{\pi\bot\ti\pi}[\exp(\varphi(|S^*|))]
=1+o(1).
\]
We next consider $a>b$. Then $\rho  =\tfrac{p\gamma}{1-p} = n^{b-a+o(1)}$.
Moreover, the condition $s^2 = \tfrac{(2-\delta)n\log n}{I}$ gives
$r=n^{2(a-b)+o(1)}$. Hence
\[
r\rho^2\le n^{o(1)},
\qquad
r\rho^4\le n^{-(a-b)+o(1)}=o(1).
\]
Propositions~\ref{prop:AFWnd-mgf} and~\ref{prop:AFWnd} yield
\[
\E_{\pi\bot\ti\pi}[\exp(\psi(\pi,\ti\pi))]=1+o(1).
\]

Let $\maP_r'=\maP_r(\,\cdot\mid\maE_\pi^1)$. In both cases, uniformly over
$r\in R_\epsilon$,
\[
\E_{\maQ}\left[\left(\frac{\diff\maP_r'}{\diff\maQ}\right)^2\right]
=\frac{\E_{\pi\bot\ti\pi}[\exp(\psi(\pi,\ti\pi))]}
{\maP_r(\maE_\pi^1)^2}
=1+o(1).
\]
Thus the conditional second moment method gives
\(
\sup_{r\in R_\epsilon}\TV(\maP_r,\maQ)=o(1).
\)
Together with~\eqref{eq:epsilon}, this yields $\TV(\maP,\maQ)=o(1)$.
\end{proof}

\begin{remark}\label{rmk:correlated-func-digraph}
    We omit a detailed review of the correlated functional digraph from~\cite{huang2024information}, which is also closely related in spirit to the lift permutation studied in~\cite{onaran2016optimal,cullina2020partial,gaudio2022exact}. The same structural reduction applies in our setting. A rigorous derivation of this reduction can be found in~\cite{huang2024information}.
\end{remark}

\subsection{Strong signal: lower bound analysis via edge density}\label{subsec:strong-signal-lwbd}

We now turn to the regime $\gamma = n^{b+o(1)}$ with a fixed constant $b>0$, under which the condition $rI = \omega(1\vee \log(1+\gamma))$ in Proposition~\ref{prop:information-density} does not hold.
Following a similar argument with~\eqref{eq:monoto-data-processing}, it suffices to consider $s^2 = \tfrac{(b\varrho^{-1}(b^{-1})-\delta)n\log n}{I}$.
In contrast to the weak-signal
regime, the quantity $\log(1+\gamma)$ is now of order $\log n$, and the main contribution to the likelihood
product comes from $\ell(1,1)$. It is therefore more natural to work with the intersection graph \begin{align}\label{eq:intersec-graph}
    H_\pi:\quad V(H_\pi)=\sfS_\pi,\quad E(H_\pi)=\bigl\{e\in \tbinom{\sfS_\pi}{2}: e\in E(G_1),\ \pi(e)\in E(G_2)\bigr\}.
\end{align}
Recall from Proposition~\ref{prop:information-density} that, under $\maP_{r}$, the edge density of $H_\pi$ is governed by the function $\varrho$. This suggests that, in the polynomial regime, the relevant truncation should be imposed on the densest subgraph of $H_\pi$, rather than on the information density considered in Section~\ref{subsec:weak-signal-lwbd}.

Fix $r\in R_\epsilon$, and let $\pi\perp\ti{\pi}$ be two independent uniform
elements of $\Pi_r$. Recall that $S^*(\pi,\ti{\pi})$ is their core set. Keep
the earlier two-sided factors $L_P$ and $L_C$, which are components of
$L_\pi L_{\ti{\pi}}$ under $\maQ$.
For the change-of-measure argument below,
define the one-sided candidate factor, for every edge path or cycle $O$, by
\begin{equation}\label{eq:one-sided-factor-preview}
  \widehat L_O
  \triangleq
  \prod_{e\in\binom{\sfS_{\ti{\pi}}}{2}\cap O}
  \ell\bigl(G_1(e),G_2(\ti{\pi}(e))\bigr).
\end{equation}
Thus
\[
  L_{\ti{\pi}}
  =\prod_{P\in\sfP^{\sfE}}\widehat L_P
   \prod_{C\in\sfC^{\sfE}}\widehat L_C.
\]
Define the set of active cycle as 
\begin{equation}\label{eq:active-orbits-preview}
  J_{\pi,\ti{\pi}}
  \triangleq
  \left\{
    C\in\sfC^{\sfE}:
    G_1(e)=G_2(\pi(e))=1
    \text{ for every }e\in C\cap\binom{\sfS_\pi}{2}
  \right\}.
\end{equation}
For any set $J$ of closed edge cycles, write $|J|\triangleq\sum_{C\in J}|C|.$ Thus $|J|$ is the total number of edges in the active cycles. On $\{C\in J_{\pi,\ti{\pi}}\}$, the one-sided factor is
$\widehat L_C=(1+\gamma)^{|C|}$. Although this can be large in the
polynomial-signal regime, the following lemma controls every inactive
cycle.

\begin{lemma}\label{lem:expect-C-notin-J}
For any fixed $\pi,\ti{\pi}\in\Pi_r$ and every closed edge cycle
$C\in\sfC^{\sfE}$,
\[
  \E_{\maP_r(\cdot\mid\pi)}
  \left[\widehat L_C\,\middle|\,
  C\notin J_{\pi,\ti{\pi}}\right]
  \le \exp(2p_{11}).
\]
\end{lemma}
\begin{proof}
   
It follows from Lemma~\ref{lem:orbit-expect} that $\E_{\maP_r(\cdot\mid\pi)}[\widehat L_C]
  =\E_{\maQ}[L_C]
  =1+\rho^{2|C|}.$
  Note that \begin{align*}
      \E_{\maP_r(\cdot\mid\pi)}[\widehat L_C] = \E_{\maP_r(\cdot\mid\pi)}[\widehat L_C\mid C\notin J_{\pi,\ti\pi}]\maP_r[C\notin J_{\pi,\ti{\pi}}\mid\pi]+\E_{\maP_r(\cdot\mid\pi)}[\widehat L_C\mid C\in J_{\pi,\ti\pi}]\maP_r[C\in J_{\pi,\ti{\pi}}\mid\pi].
  \end{align*}
Moreover,
\[
  \maP_r[C\in J_{\pi,\ti{\pi}}\mid\pi]=p_{11}^{|C|},
  \qquad
  \widehat L_C=(1+\gamma)^{|C|}
  \quad\text{on }\{C\in J_{\pi,\ti{\pi}}\}.
\]
Consequently,
\begin{align*}
  &\E_{\maP_r(\cdot\mid\pi)}
  \left[\widehat L_C\,\middle|\,
  C\notin J_{\pi,\ti{\pi}}\right]=
  \frac{1+\rho^{2|C|}-p_{11}^{|C|}(1+\gamma)^{|C|}}
       {1-p_{11}^{|C|}}
  \le \frac{1}{1-p_{11}^{|C|}}
  \le \exp(2p_{11}^{|C|})
  \le \exp(2p_{11}).
\end{align*}
The first inequality follows from
$ p_{11}(1+\gamma)
  =\bigl(p+\rho-\rho p\bigr)^2\ge\rho^2$
and the second holds because $p_{11}\le1/2$.
\end{proof}

We next establish the corresponding product bound. Fix a constant
$\kappa_b = 2\varrho^{-1}(b^{-1})$ such that $rp_{11}\le\kappa_b$ for all sufficiently large
$n$.

\begin{lemma}\label{lem:large-gamma-ding-du}
For any $\maE_\pi$ such that
$\maP_r(\maE_\pi\mid\pi)>0$,,
    \begin{align*}
        \E_{\pi\bot\ti\pi} \E_{\maQ}[L_\pi L_{\ti\pi} \indc{\maE_\pi}\indc{\maE_{\ti\pi}}]\le \E_{\pi\bot\ti\pi}\bigg[\E_{J_{\pi,\ti\pi}\sim\maP_{r}(\cdot \mid \pi,\maE_\pi)}\Big[ \frac{e^{\kappa_b|S^*|}(1+\gamma)^{|J_{\pi,\ti\pi}|}}{\maP_{r}(\maE_\pi\mid \pi,J_{\pi,\ti\pi})}\Big]\bigg].
    \end{align*}
\end{lemma}
\begin{proof}
   For fixed $(\pi,\ti{\pi})$, change measure with respect to the true
alignment $\pi$ yields
\begin{align*}
  &\E_{\maQ}\left[
    L_\pi L_{\ti{\pi}}
    \indc{\maE_\pi}\indc{\maE_{\ti{\pi}}}
  \right]=
  \E_{\maP_r(\cdot\mid\pi)}\left[
    L_{\ti{\pi}}
    \indc{\maE_\pi}\indc{\maE_{\ti{\pi}}}
  \right]
  \le
  \E_{\maP_r(\cdot\mid\pi)}\left[
    L_{\ti{\pi}}\indc{\maE_\pi}
  \right].
\end{align*}
Conditional on $J_{\pi,\ti{\pi}}=J$, decompose the one-sided likelihood
$L_{\ti{\pi}}$. The active cycles contribute $(1+\gamma)^{|J|}$, every
path factor has mean one, and Lemma~\ref{lem:expect-C-notin-J} bounds each inactive cycle factor.
Therefore,
\begin{align*}
  \E_{\maP_r(\cdot\mid\pi)}
  [L_{\ti{\pi}}\mid J_{\pi,\ti{\pi}}=J]
  &\le
  (1+\gamma)^{|J|}
  \exp\left(2p_{11}\binom{|S^*|}{2}\right)\le
  (1+\gamma)^{|J|}e^{\kappa_b|S^*|},
\end{align*}
where the first inequality holds because there are at most $\binom{|S^*|}{2}$ inactive cycles.
Finally, Bayes' rule gives
\[
  \maP_r(J\mid\pi)
  =
  \frac{
    \maP_r(\maE_\pi\mid\pi)
    \maP_r(J\mid\pi,\maE_\pi)
  }{
    \maP_r(\maE_\pi\mid\pi,J)
  }.
\]
Substituting the preceding bound and averaging over $(\pi,\ti{\pi})$ proves
the result.
\end{proof}

We now introduce the event $\maE_\pi^2$. Since $s^2 \le (b\varrho^{-1}(b^{-1})-\delta)n\log n/I$ and $I=(b+o(1))p_{11}\log n$, we may choose $\epsilon>0$ sufficiently small such that $r p_{11}\le \varrho^{-1}(b^{-1})-\delta/2$ uniformly for all $r\in R_\epsilon$. Choose a small constant $\xi$ such that $\varrho(\varrho^{-1}(b^{-1})-\delta/2)<\xi<b^{-1}$, and then choose a small constant $\eta>0$ such that $10\eta<1-b\xi$.
For any graph $H$ with $r$ vertices, we define $H\in \maH(r)$ if  1) \begin{align}\label{eq:ratio-xi}
    \max_{\emptyset\neq V\subseteq V(H)} \frac{e_H(V)}{|V|}\le \xi;
\end{align} 2) the maximum degree of $H$, defined as
$\Delta(H)\triangleq \max_{v\in V(H)}
|\{u\in V(H):uv\in E(H)\}|$, satisfies
$\Delta(H)\le \log n$; 3) for every $k\ge 3$, the number of labeled $k$-cycles in $H$ is at most $n^{\eta k}$.

For any $\pi\in\Pi_r$, let $H_\pi$ denote the graph with vertex set $V(H_\pi)=\sfS_\pi$ and edge set $E(H_\pi)=\{e:e\in E(G_1),\,\pi(e)\in E(G_2)\}$. Let $\maE_\pi^2$ denote the event that $H_\pi\in \maH(r)$. Since $rI/(b\log n)<\varrho^{-1}(\xi)$, Proposition~\ref{prop:information-density} implies that the first condition holds with probability $1-o(1)$ under $\maP_{r}$ for $H_\pi$.
For the maximum degree, each degree in $H_\pi$ is distributed as $\Bin(r-1,p_{11})$, and $rp_{11}=O(1)$. Hence, by a standard Chernoff bound (see, e.g.,~\cite[Theorem 4.4]{mitzenmacher_probability_2005}) and a union bound,
\begin{align*}
    \maP_{r}\big(\Delta(H_\pi)>\log n\big)
    &\le r\,\prob{\Bin(r-1,p_{11})>\log n}  \le r\left(\frac{e(r-1)p_{11}}{\log n}\right)^{\log n}
     = o(1).
\end{align*}
For the cycle count, fix $k\ge 3$. The expected number of labeled $k$-cycles in $H_\pi$ is at most $(r)_k p_{11}^k\le (rp_{11})^k=O(1)^k$. Therefore, by Markov's inequality,
\[
\P\left[\text{the number of labeled $k$-cycles in }H_\pi>n^{\eta k}\right]
\le O(1)^k n^{-\eta k}.
\]
Summing over $3\le k\le r$ gives $\sum_{k=3}^r O(1)^k n^{-\eta k}=o(1)$. Combining the above estimates, we obtain $\maP_{r}(\maE_\pi^2)=1-o(1)$. 
{For the active edge-cycle set $J_{\pi,\ti\pi}$, define
\[
E(J_{\pi,\ti\pi})
\triangleq
\bigcup_{C\in J_{\pi,\ti\pi}}
\left(C\cap\binom{\sfS_\pi}{2}\right),
\qquad
V(J_{\pi,\ti\pi})
\triangleq
\bigcup_{e\in E(J_{\pi,\ti\pi})}V(e),
\]
and let
\[
H(J_{\pi,\ti\pi})
\triangleq
\bigl(V(J_{\pi,\ti\pi}),E(J_{\pi,\ti\pi})\bigr).
\]
By
construction, it is invariant under
$\theta=\ti\pi^{-1}\circ\pi$, satisfies
$|E(J_{\pi,\ti\pi})|=|J_{\pi,\ti\pi}|$, and is a subgraph of $H_\pi$.
In particular, on $\maE_\pi^2$, we have  $H(J_{\pi,\ti\pi})\in \maH(|V(J_{\pi,\ti\pi})|)$.}
We are now ready to establish the main result in this section.

\begin{proposition}\label{prop:main-IT-lwbd-poly-gamma}
Given $\pi\in \Pi_r$ with $H_\pi\in \maH(r)$, then there exists a small constant such that, for any $\ti\pi\in \Pi_r$ and $J_{\pi,\ti\pi}\sim \maP_r(\cdot \mid \pi,\maE_\pi^2)$, we have $\maP_{r}(\maE_\pi^2\mid\pi,J_{\pi,\ti\pi})\ge (1-o(1))c^{|J_{\pi,\ti\pi}|}$. Moreover, when $s^2\le \tfrac{(b\varrho^{-1}(b^{-1})-\delta)n\log n}{I}$, we have \begin{align*}
        \E_{\pi\bot\ti\pi} \Big[e^{\kappa_b|S^* |} [c^{-1}(1+\gamma)]^{|E(J_{\pi,\ti\pi})|}\Big] = 1+o(1).
    \end{align*}
\end{proposition}

The proof of Proposition~\ref{prop:main-IT-lwbd-poly-gamma} is deferred to Appendix~\ref{subsec:proof-prop-main-IT-lwbd}. Combining Lemma~\ref{lem:large-gamma-ding-du} with Proposition~\ref{prop:main-IT-lwbd-poly-gamma} gives  $\TV(\maP_r,\maQ) = o(1)$ when $s^2\le \tfrac{(b\varrho^{-1}(b^{-1})-\delta)n\log n}{I}$. Hence, we prove the lower bound when $\gamma = n^{b+o(1)}$ in Theorem~\ref{thm:main-IT}.

\begin{remark}
The extra factor $e^{\kappa_b|S^*|}$ in Lemma~\ref{lem:large-gamma-ding-du} is one of the main differences between our argument and that of~\cite{ding2023detection}. In their model, the analogue of the inactive-orbit contribution is bounded by $1$, using the assumptions $1+\gamma\ge 10$ and $p(1+\gamma)\le 0.1$, which are compatible with their setting. In our subgraph sampling model, however, these assumptions are substantially more restrictive than the parameter range considered here. In particular, they are not implied by the polynomial-signal assumption $0<b\le a$; near the endpoint $a=b$, the quantity $p(1+\gamma)$ need not be small. Thus, inactive cyclic orbits in our setting cannot be uniformly bounded by $1$. Instead, we use the softer single-cycle estimate $\exp(2p_{11})$. After summing over inactive cyclic orbits, this yields the factor $e^{\kappa_b|S^*|}$ in Lemma~\ref{lem:large-gamma-ding-du}. This additional term is then absorbed into the counting argument and does not affect the final threshold.
\end{remark}

\section{Information-theoretic upper bounds}\label{sec:upbd}


\subsection{Matching via log-likelihood statistic for subpolynomial signals}

In this section, we consider the subpolynomial-signal regime where $\gamma\le n^{o(1)}$. Let $\Phi_\epsilon$ be the collection of all injective mappings $\phi:S\subseteq V(G_1)\mapsto T\subseteq V(G_2)$ such that $|S|=|T|=k\triangleq (1-\epsilon){s^2}/{n}.$
We consider the maximal log-likelihood statistic
\[
\maT_1 = \max_{\phi\in \Phi_\epsilon} Z_\phi^\ell\triangleq \max_{\phi\in\Phi_\epsilon}\sum_{\{u,v\}\subseteq S}
\log \ell\big(G_1(uv),G_2(\phi(u)\phi(v))\big).
\]
Under $\maP$, since $\prob{|\sfS_\pi|\ge k}=1-o(1)$, with high probability there exist subsets $S^*\subseteq \sfS_\pi$ and $T^*=\pi(S^*)\subseteq \sfT_\pi$ with $|S^*|=|T^*|=k$. We then take the true injective mapping $\phi^*=\pi|_{S^*}:S^*\mapsto T^*$ as the signal alignment. Under $\maQ$, we control the statistic by a union bound over $\Phi_\epsilon$, whose cardinality is $|\Phi_\epsilon|=\binom{s}{k}^2 k!.$
Specifically, we have the following proposition.
\begin{proposition}\label{prop:upbd-weak-signal}
    When $\gamma \le  n^{o(1)}$ and $s^2 \ge \tfrac{(2+\delta)n\log n}{I}$, we have $\TV(\maP,\maQ) = 1-o(1)$.
\end{proposition}
\begin{proof}
    Pick $\tau = (1+\delta/4)k\log n$. It suffices to show $\maQ(\maT_1\ge \tau)+\maP(\maT_1<\tau) = o(1)$. Under $\maQ$, we have $\expect{\exp(Z_\phi^\ell)} = 1$ for any $\phi\in\Phi_\epsilon$. Consequently, \begin{align*}
        \maQ(Z_\phi^\ell\ge (1+\delta/4)k\log n) = \maQ(\exp(Z_\phi^\ell)\ge n^{(1+\delta/4)k})\le n^{-(1+\delta /4)k}.
    \end{align*}
    Applying union bound yields \begin{align*}
        \maQ(\maT_1\ge \tau)\le |\Phi_\epsilon| n^{-(1+\delta/4)k} =\binom{s}{k}^2 k!  n^{-(1+\delta/4)k}\le (e^2s^2/k)^kn^{-(1+\delta/4)k}\le n^{-\delta k/8} = o(1).
    \end{align*}
    Under $\maP$, we have $\expect{Z_{\phi^*}^\ell} = \binom{k}{2}I\ge (1+\delta /3)k\log n$. By Lemma~\ref{lem:var-gamma-bounded}, we have $$\var[Z_{\phi^*}^\ell]\le (2+\log (1+\gamma))\binom{k}{2}I\le (2+\log(1+\gamma)) k^2,$$
    where the last inequality holds because $G_1(e),G_2(\pi(e))\in \sth{0,1}$ and the mutual information $I= I(G_1(e);G_2(\pi(e)))\le \log 2$. By Chebyshev's inequality,\begin{align*}
        \maP(Z_{\phi^*}^\ell<\tau)\le \maP\Big(Z_{\phi^*}^\ell-\expect{Z_{\phi^*}^\ell}<-\frac{\delta k\log n}{12}\Big)\le \frac{(2+\log(1+\gamma)) k^2}{\delta^2 k^2 \log^2n/144} = o(1),
    \end{align*}where the last equality follows from $\log(1+\gamma) = o(\log n)$. This implies $\maP(\maT_1<\tau)=o(1)$, which completes the proof.
\end{proof}

\subsection{Matching via densest subgraph for polynomial signals}

We next prove the information-theoretic upper bound in the polynomial regime. In this case the likelihood ratio is dominated by $\ell(1,1)$, so the statistic from previous subsection is no longer the most natural choice. Instead, we scan over partial alignments of size comparable to the overlap and test whether the corresponding common-edge graph contains an unusually dense subgraph. 
Recall that $k = (1-\epsilon)s^2/n$ and $\Phi_\epsilon$ denote the collection of all injective mappings $\phi:S\subseteq V(G_1)\mapsto T\subseteq V(G_2)$ with $|S| = |T| = k$.
We consider the densest subgraph statistic\begin{align*}
    \maT_2 = \max_{\phi\in \Phi_\epsilon}\max_{\emptyset \neq V\subseteq S} \frac{Z_\phi(V)}{|V|}  = \max_{\phi\in \Phi_\epsilon}\max_{\emptyset \neq V\subseteq S} \frac{|\{e\in \binom{V}{2}:e\in E(G_1),\phi(e)\in E(G_2)\} |}{|V|}.
\end{align*}
Specifically, we have the following proposition.
\begin{proposition}\label{prop:upbd-strong-signal}
    When $\gamma = n^{b+o(1)}$ for some constant $b>0$ and $s^2\ge \tfrac{(b\varrho^{-1}(b^{-1})+\delta)n\log n}{I}$, we have $\TV(\maP,\maQ) = 1-o(1)$.
\end{proposition}
\begin{proof}
    By monotonicity, we only need to consider  $s^2= \tfrac{(b\varrho^{-1}(b^{-1})+\delta)n\log n}{I}$.
    Pick a constant $\tau$ such that $b^{-1}<\tau<\varrho(\varrho^{-1}(b^{-1})+\delta/2b)$. Under $\maP$, since $\prob{|\sfS_\pi|\ge k}=1-o(1)$, with high probability there exist subsets $S^*\subseteq \sfS_\pi$ and $T^*=\pi(S^*)\subseteq \sfT_\pi$ with $|S^*|=|T^*|=k$. Let $\phi^* = \pi\mid_{S^*}:S^*\mapsto \pi(S^*)$. Following the notation in~\eqref{eq:intersec-graph}, the intersection graph $H_{\phi^*}\sim \maG(k,p_{11})$. Since $I = (b+o(1))p_{11}\log n$, we have \begin{align*}
        k p_{11}\ge (1-\epsilon)  \frac{s^2 p_{11}}{n}\ge (1-\epsilon)  \frac{b\varrho^{-1}(b^{-1})+\delta}{b+o(1)}\ge \varrho^{-1}(b^{-1})+\frac{\delta}{2b},
    \end{align*}where the last inequality holds when we pick constant $\epsilon$ sufficiently small. It follows from Proposition~\ref{prop:information-density} that $\maP(\maT_2<\tau) = o(1)$.

    We then show $\maQ(\maT_2\ge\tau) = o(1)$. Under $\maQ$, let $m = |V|$. Then $Z_\phi(V)\sim \Bin(\binom{m}{2},p^2)$. By Chernoff bound (see, e.g.,~\cite[Theorem 4.4]{mitzenmacher_probability_2005}),\begin{align*}
        \maQ(Z_\phi(V)\ge \tau m)\le \exp\pth{\tau m\log \frac{e\binom{m}{2}p^2}{\tau m}}\le \exp\Big(\tau m\log \frac{ep^2 m}{2\tau}\Big).
    \end{align*}
Since there are at most $\binom{s}{m}^2 m!$ injective mappings $\phi$ with domain and image of size $m$, a union bound gives
    \begin{align*}
        \maQ(\maT_2\ge \tau)&\le \sum_{m\ge 1}\binom{s}{m}^2 m!\exp\Big(\tau m\log \frac{ep^2 m}{2\tau}\Big)\\&\overset{\mathrm{(a)}}{\le}\sum_{m\ge 1} \exp\pth{m[(\tau-1)\log m+(1+2a-b-2\tau a)\log n+o(\log n)]}\\&\overset{\mathrm{(b)}}{\le}\sum_{m\ge 1} \exp(m(1-\tau b+o(1))\log n) = o(1),
    \end{align*}
    where $\mathrm{(a)}$ is because $\binom{s}{m}\le (\tfrac{es}{m})^{m}$, $s^2 = n^{1+2a-b+o(1)}$, and $p = n^{-a+o(1)}$; $\mathrm{(b)}$ uses $m\le k\le n^{2a-b+o(1)}$; the last equality follows from  constant $1-\tau b<0$. Hence $\maQ(\maT_2\ge \tau)+\maP(\maT_2<\tau) = o(1)$ and $\TV(\maP,\maQ) =1- o(1)$.
\end{proof}

\subsection{Phase transition from strong signal to weak signal}

We now finish the proof of the sharp information-theoretic thresholds established in Theorem~\ref{thm:main-IT}. In view of Theorem~\ref{thm:main-IT}, the sharp constant changes from $2$ in the case $b=0$ to $b\varrho^{-1}(b^{-1})$ in the polynomial regime $b>0$. A natural question is whether this transition is continuous as $b\to 0_+$. In this subsection, we show that the transition is indeed smooth, namely $\lim_{b\to 0_+} b\varrho^{-1}(b^{-1})=2$. Specifically, we prove the following proposition.
\begin{proposition}\label{prop:phase-transition}
For the densest subgraph function defined in~\eqref{eq:densest-func}, we have $\lim_{b\to 0_+} b\varrho^{-1}(b^{-1}) = 2$.
\end{proposition}

\begin{proof}
    
    We first prove $\liminf_{\lambda\to \infty} \frac{\varrho(\lambda)}{\lambda} \ge \frac{1}{2}.$
    For \(H\sim\mathcal G(N,\tfrac{\lambda}{N})\), we have \(|E(H)|\sim\operatorname{Bin}\bigl(\binom{N}{2},\tfrac{\lambda}{N}\bigr)\), and, for each fixed \(\lambda\),
\[
\mathbb E\!\left[\frac{|E(H)|}{N}\right]
=\frac{\lambda}{2}\left(1-\frac1N\right),
\qquad
\operatorname{Var}\!\left(\frac{|E(H)|}{N}\right)
=\frac{\lambda(N-1)}{2N^2}\left(1-\frac{\lambda}{N}\right)
=O_\lambda(N^{-1}).
\]Hence, by Chebyshev’s inequality, \(|E(H)|/N\to\lambda/2\) in probability, and therefore, for every fixed \(\epsilon>0\),
    \begin{align*}
        \prob{\frac{|E(H)|}{N}>(1-\epsilon)\frac{\lambda}{2}} = 1-o(1)\quad \text{  as  } N\to\infty,
    \end{align*}which implies $\varrho(\lambda)>(1-\epsilon)\lambda /2$. Letting $\epsilon\to 0_+$ gives $\liminf_{\lambda\to \infty} \frac{\varrho(\lambda)}{\lambda} \ge \frac{1}{2}.$

    We then prove $\limsup_{\lambda\to \infty} \frac{\varrho(\lambda)}{\lambda} \le \frac{1}{2}.$  For any $V\subseteq V(H)$ with $|V| = k$, we have $|E(H[V])|\sim \Bin(\binom{k}{2},\frac{\lambda}{N})$. For any constant $\epsilon>0$, by Chernoff bound,\begin{align*}
        \prob{|E(H[V])|\ge (1+\epsilon) \frac{\lambda k}{2}}&\le \exp\pth{-(1+\epsilon)\frac{\lambda k}{2}\log \frac{\frac{(1+\epsilon)\lambda k}{2}}{\binom{k}{2}\frac{\lambda }{N}}+(1+\epsilon)\frac{\lambda k}{2} -\binom{k}{2}\frac{\lambda }{N} }\\&=\exp\pth{-\frac{\lambda k}{2} f_\epsilon\Big(\frac{k-1}{N}\Big)},
    \end{align*}where $f_\epsilon(x) = (1+\epsilon)\log((1+\epsilon)/x)-(1+\epsilon)+x$. Since $f_\epsilon'(x) = 1-\tfrac{1+\epsilon}{x}<0$ for $0<x\le 1$, we have $f_\epsilon(x)\ge f_\epsilon(1)>0$. Since $\log (e/x)\ge 0$ when $0<x\le 1$ and $\tfrac{f_\epsilon(x)}{\log (e/x)}\to 1+\epsilon$ as $x\to 0_+$, we have $\tfrac{f_\epsilon(x)}{\log (e/x)}\ge c_\epsilon$ for all $0<x\le 1$ with some constant $c_\epsilon>0$. Then, when $\lambda>\tfrac{8}{c_\epsilon}$, we have \begin{align*}
        &~\prob{\exists V\subseteq V(H): |E(H[V])| \ge \frac{(1+\epsilon)\lambda |V|}{2}}\le \sum_{k= 2}^N\binom{N}{k}\exp\pth{-\frac{\lambda k}{2} f_\epsilon\Big(\frac{k-1}{N}\Big)}\\\le&~\sum_{k= 2}^N \exp\pth{k\log (eN/k)-4k\log (eN/(k-1))}\le N(1/e N)^3 = o(1),
    \end{align*}where the last inequality follows from the fact that $k\mapsto k\log(k/eN)$ has derivative $\log(k/N)<0$. Letting $\epsilon\to 0_+$ gives $\limsup_{\lambda\to \infty} \frac{\varrho(\lambda)}{\lambda} \le \frac{1}{2}.$ Combining this with $\liminf_{\lambda\to \infty} \frac{\varrho(\lambda)}{\lambda} \ge \frac{1}{2}$ yields \begin{align}\label{eq:0.5-limit}
        \lim_{\lambda\to\infty} \frac{\varrho(\lambda)}{\lambda} = \frac{1}{2}
    \end{align} which also gives $\lim_{\lambda\to\infty} = \varrho(\lambda) = \infty$.
    Recall the densest subgraph function $\varrho$ defined in~\eqref{eq:densest-subgraph} has the inverse function $\varrho^{-1}:[1,\infty)\to[1,\infty)$. We then have $\lim_{b\to 0} \varrho^{-1}(b^{-1}) \to \infty$. Combining this with~\eqref{eq:0.5-limit} yields $$\lim_{b\to 0}\frac{\varrho(\varrho^{-1}(b^{-1}))}{\varrho^{-1}(b^{-1})} = \frac{1}{2}.$$ This is equivalent to $\lim_{b\to 0_+} b\varrho^{-1}(b^{-1}) = 2$. We prove the proposition.
\end{proof}







\section{Discussion}

 \emph{Sharp sample complexity and the core set.}
We study correlation detection when only two randomly sampled induced
subgraphs are observed. Theorem~\ref{thm:main-IT} gives the sharp
leading-order information-theoretic threshold
\[
    s_\star^2
    = C_{\varrho}(b)\frac{n\log n}{I}.
\]
The core set $S^*(\pi,\widetilde{\pi})$ plays a central role in the
lower-bound analysis. For two candidate partial alignments, it collects
the vertex cycles on which the two alignments agree. All closed edge
orbits in the likelihood-product decomposition are generated by the
induced action on $\binom{S^*}{2}$. The second-moment calculation can
therefore be reduced to controlling the size and cycle structure of the
core set. This reduction allows us to retain the sharp threshold
constant.
  
 \emph{A general conditional second-moment framework.}
We also prove a general result for the conditional second-moment
method. If two probability measures are asymptotically indistinguishable
in total variation, then there exists a high-probability event on which
the conditional likelihood-ratio second moment approaches its null
value. The converse follows from the standard second-moment bound.
Thus, at the level of existence, conditioning causes no loss of
information. In a specific model, one still needs to find a tractable
conditioning event and control the resulting second moment.
 
 \emph{Strong signals and densest subgraphs.}
The behavior changes in the polynomial strong-signal regime
$\gamma=n^{b+o(1)}$ with $b>0$. In this regime, common present edges
dominate the likelihood. The relevant object is therefore the
intersection graph induced by a candidate alignment. Its maximum
edge--vertex ratio is described by the sparse Erd\H{o}s--R\'enyi
densest-subgraph function $\varrho$. This gives the sharp constant
\[
    C_{\varrho}(b)=b\varrho^{-1}(b^{-1}).
\]
The same quantity appears in both the upper- and lower-bound arguments.
Proposition~\ref{prop:phase-transition} further shows that
\[
    \lim_{b\downarrow0}b\varrho^{-1}(b^{-1})=2.
\]
Hence, although the proof mechanism changes at $b=0$, the threshold
constant remains continuous. Equivalently,
$\varrho(\lambda)/\lambda\to1/2$ as $\lambda\to\infty$.

\section*{Acknowledgment}
We warmly thank Shuyang Gong, Zhangsong Li and Jianqiao Wang for discussions and their
valuable advice. An earlier version of this work established rate-optimal sample-complexity bounds in most parameter regimes, generally without identifying the sharp leading constants. In preparing the present version, the authors used OpenAI's GPT-5.4 Pro, GPT-5.5 Pro, and GPT-5.6 Pro as research-assistance tools to explore, check, and refine the arguments leading to the improved lower bound constants in Sections~\ref{subsec:weak-signal-lwbd} and~\ref{subsec:strong-signal-lwbd}. With the exception of this AI-assisted refinement, the mathematical ideas, analyses, and proofs in the paper were developed independently by the authors. These models were also used for language editing and stylistic polishing. All AI-assisted suggestions and derivations were independently examined, verified, and revised by the authors, who take full responsibility for the correctness and content of the entire manuscript.

\appendix

\section{Proof of Propositions}

\subsection{Proof of Proposition~\ref{prop:AFWnd}}\label{apd:proof-prop-AFWnd}
For any $k\ge 1$, define $\sfC_k^\sfE\triangleq \sth{C\in \sfC^\sfE:|C| = k}$ as the set of cycle with length $k$. 
Let
\[
W_{\mathrm{nd}}(\pi,\ti\pi)
\triangleq
\prod_{C\in\sfC^\sfE_{\mathrm{nd}}}(1+\rho^{2|C|}) = (1+\rho^2)^{|\sfC_2^\sfE|-\binom{C_1}{2}}\prod_{k\ge 2}(1+\rho^{2k})^{|\sfC_k^\sfE|},
\]
where $\sfC^\sfE_{\mathrm{nd}}\triangleq (\cup_{k\ge 1}\sfC_k^\sfE)\backslash \binom{\sfC^\sfV_1}{2}$ contains all closed edge cycles except
the one-edge cycles between two fixed points. Here $\sfC_1^\sfV$ denotes the set of fixed points.
We then have \begin{align}\label{eq:fixed-orbit}
\exp(\psi(\pi,\ti\pi))
\le \exp(\varphi(C_1))W_{\mathrm{nd}}(\pi,\ti\pi).
\end{align}

The case $\rho=0$ is immediate, so assume $\rho>0$. The condition
$r\rho^4=o(1)$ implies $\rho=o(1)$. Split the non-fixed closed edge cycles
into three classes and write $W_{\mathrm{nd}}=W_1W_2W_3$. Here $W_1$
contains the cycles internal to one  vertex cycle, $W_2$ those
joining a fixed point to a  vertex cycle, and $W_3$ those joining
two distinct vertex cycles. Since \(W_i\ge 1\) for each \(i\), the AM--GM inequality applied to \(W_1^3,W_2^3,W_3^3\) gives
\begin{align}\label{eq:Wnd}
W_{\mathrm{nd}}-1\le\frac13\sum_{i=1}^3(W_i^3-1).
\end{align}
It is therefore enough to show that, for $i=1,2,3$,
\begin{align}\label{eq:Wi}
\E\left[\exp(\varphi(C_1))(W_i^3-1)\right]=o(1).
\end{align}

\emph{The $W_1$ term.}
A vertex cycle of length $k\ge2$ has at most $k$ internal unordered-edge
orbits, each of length at least $\lceil k/2\rceil$. Hence
\[
W_1^3\le
\exp\Big(\sum_{k\ge2}3k\rho^{2\lceil k/2\rceil}C_k\Big).
\]
Since $\rho=o(1)$, we have
$3k\rho^{2\lceil k/2\rceil}\le 1$. Lemma~\ref{lem:psi} with $u_\ell = 3\ell \rho^{2\lceil k/2\rceil}$ therefore gives $$\expect{\indc{C_1 =k} (W_1^3-1)}\le \frac{(r^2/s^2)^{k}}{k!}(\exp(\Psi_1)-1),$$ where 
\[
\Psi_1\triangleq
\sum_{\ell\ge2}\frac{(r^2/s^2)^\ell}{\ell}
\left(\exp(3\ell\rho^{2\lceil \ell/2\rceil})-1\right)
\le6\sum_{\ell\ge2}(r^2/s^2)^\ell\rho^{2\lceil \ell/2\rceil}
\le C(r^2/s^2)^2\rho^2=o(1).
\]
Under the weak signal regime $\gamma \le n^{o(1)}$, it follows from~\eqref{eq:asy-I} that $I = n^{-2(a-b)+o(1)}$ and $\rho = n^{b-a+o(1)}$. Thus $r \rho^2 \le \tfrac{2(1+\epsilon)\log n}{I}\rho^2\le n^{o(1)}$.
Recall from the two estimates in the proof of Proposition~\ref{prop:AFWnd-mgf}, for
$2\le k\le1/\log(1+\rho^2)$ and $k>1/\log(1+\rho^2)$, respectively, that
\[
\sum_{k=0}^{r}\exp(\varphi(k))\frac{(r^2/s^2)^k}{k!}=O(1).
\]
Therefore
\begin{align}\label{eq:W1}
\E\left[\exp(\varphi(C_1))(W_1^3-1)\right]
\le(\exp(\Psi_1)-1)
\sum_{k=0}^{r}\exp(\varphi(k))\frac{(r^2/s^2)^k}{k!}=o(1),
\end{align}
which proves~\eqref{eq:Wi} for $i=1$.

\emph{The $W_2$ term.}
A fixed point and a vertex cycle of length $\ell\ge2$ generate one edge
orbit of length $\ell$. Thus
\[
W_2^3\le
\exp\Big(3C_1\sum_{\ell\ge2}\rho^{2\ell}C_\ell\Big).
\]
For $C_1=k$, Lemma~\ref{lem:psi} with $u_\ell = 3\rho^{2\ell}$ gives
\[
\E\left[\indc{C_1=k}(W_2^3-1)\right]
\le\frac{(r^2/s^2)^k}{k!}(\exp(\Psi_2(k))-1),
\]
where
\[
\Psi_2(k)\triangleq
\sum_{\ell\ge2}\frac{(r^2/s^2)^\ell}{\ell}
\left(\exp(3k\rho^{2\ell})-1\right).
\]
Uniformly for $0\le k\le r$,
\[
\Psi_2(k)\le3k\exp(3r\rho^4)
\sum_{\ell\ge2}\frac{((r^2/s^2)\rho^2)^\ell}{\ell}
\le Cr(r^2/s^2)^2\rho^4=o(1).
\]
Here the last step uses $e^{3r\rho^4}=1+o(1)$,
$r\rho^4=o(1)$, and $r^2/s^2\le1$.
Consequently,
\begin{align}\label{eq:W2}
\expect{\exp(\varphi(C_1))(W_2^3-1)}
&\le \sup_{0\le k\le r}(\exp(\Psi_2(k))-1)
\sum_{k=0}^{r}\exp(\varphi(k))\frac{(r^2/s^2)^k}{k!}
=o(1).
\end{align}

\emph{The $W_3$ term.}
Two vertex cycles of lengths $k$ and $\ell$ generate
$g=\gcd(k,\ell)$ edge orbits, each of length
$L=\operatorname{lcm}(k,\ell)$. Since
$3g\rho^{2L}\le3r\rho^4=o(1)$, for all sufficiently large $n$,
\begin{align}\label{eq:cross}
(1+\rho^{2L})^{3g}-1
\le6\min\{k,\ell\}\rho^{2\max\{k,\ell\}}.
\end{align}
For cycles of different lengths, the sum of the right-hand side of
\eqref{eq:cross} is at most
\[
6\sum_{2\le k<\ell}kC_kC_\ell\rho^{2\ell}
\le6r\sum_{\ell\ge2}C_\ell\rho^{2\ell},
\]
where the inequality uses $\sum_{k\ge 1} kC_k = |S^*|\le r$ from~\eqref{eq:core-vertex-cycle-decomposition}.
For two distinct cycles of the same length, it is at most
\[
6\sum_{\ell\ge2}\ell\binom{C_\ell}{2}\rho^{2\ell}
\le3r\sum_{\ell\ge2}C_\ell\rho^{2\ell}.
\]
Thus, using $\log(1+z)\le z$,
\[
W_3^3\le
\exp\Big(9r\sum_{\ell\ge2}\rho^{2\ell}C_\ell\Big).
\]
Lemma~\ref{lem:psi} with $u_\ell=9r\rho^{2\ell}$ gives
\begin{align*}
    \expect{\indc{C_1 =k} (W_3^3-1)}\le \frac{(r^2/s^2)^{k}}{k!}(\exp(\Psi_3)-1),
\end{align*}
where
\[
\Psi_3
\triangleq
\sum_{\ell\ge2}\frac{(r^2/s^2)^\ell}{\ell}
\left(\exp(9r\rho^{2\ell})-1\right)
\le9r\exp(9r\rho^4)
\sum_{\ell\ge2}\frac{((r^2/s^2)\rho^2)^\ell}{\ell}
\le Cr(r^2/s^2)^2\rho^4=o(1).
\]
Again, the last step uses $e^{9r\rho^4}=1+o(1)$,
$r\rho^4=o(1)$, and $r^2/s^2\le1$.
It follows that
\begin{align}\label{eq:W3}
\E\left[\exp(\varphi(C_1))(W_3^3-1)\right]
&\le(\exp(\Psi_3)-1)
\sum_{k=0}^{r}\exp(\varphi(k))\frac{(r^2/s^2)^k}{k!}
=o(1).
\end{align}
Combining~\eqref{eq:fixed-orbit},
\eqref{eq:Wnd}, and~\eqref{eq:W1}--\eqref{eq:W3}
proves Proposition~\ref{prop:AFWnd}.

\subsection{Proof of Proposition~\ref{prop:main-IT-lwbd-poly-gamma}}\label{subsec:proof-prop-main-IT-lwbd}

In this subsection, we prove Proposition~\ref{prop:main-IT-lwbd-poly-gamma}.
We first introduce the following lemma on $\maP_r(\maE_\pi^2\mid \pi,J_{\pi,\ti\pi})$.

\begin{lemma}\label{lem:conditional-prob-lwbd}Fix $\pi,\ti\pi$ and a realization $J$ in the support of
$J_{\pi,\ti\pi}$ under
$\maP_r(\cdot\mid\pi,\maE_\pi^2)$. When $rp_{11} = O(1)$, there exists a small constant $c>0$ such that, 
    $\maP_{r}(\maE_\pi^2\mid \pi,J_{\pi,\ti\pi} = J)\ge (1-o(1))c^{|E(J)|}$.
\end{lemma}
The proof of Lemma~\ref{lem:conditional-prob-lwbd} is deferred to Section~\ref{subsec:proof-lem-condition-prob}.
Then, it suffices to show $\E_{\pi\bot\ti\pi} \Big[e^{\kappa_b|S^* |} [c^{-1}(1+\gamma)]^{|E(J_{\pi,\ti\pi})|}\Big] = 1+o(1)$ on the information-theoretic lower bound for $\gamma = n^{b+o(1)}$ regime in Theorem~\ref{thm:main-IT}. 
Let $H_\pi$ denote the graph with vertex set $V(H_\pi)=\sfS_\pi$, where an edge
$e$ belongs to $E(H_\pi)$ if and only if $e\in E(G_1)$ and $\pi(e)\in E(G_2)$.
Since $\pi\perp \ti\pi$ and both are uniformly distributed over $\Pi_r$, by symmetry we may fix $\pi$ and take $\ti\pi\sim \mathrm{Unif}(\Pi_r)$. 
The following proposition gives an upper bound on
\(
\E_{\pi\perp\ti\pi}\Big[
    e^{\kappa_b|S^*|}
    \bigl(c^{-1}(1+\gamma)\bigr)^{|E(J_{\pi,\ti\pi})|}
\Big],
\)
and its proof is deferred to Section~\ref{subsec:proof-prop-IT-lwbd-mid}.
\begin{proposition}\label{prop:prop-IT-lwbd-mid} 
For any $\pi\in \Pi_r$ with $H_\pi\in \maH(r)$, we have 
    \begin{align}\label{eq:lemma-str-snt}
        \frac{1}{|\Pi_r|}\sum_{\ti\pi\in \Pi_r} e^{\kappa_b|S^* |} [c^{-1}(1+\gamma)]^{|E(J_{\pi,\ti\pi})|}\le \frac{1}{1-e^{\kappa_b}r^2/s^2} \exp(2S_{\mathrm{tr}}(H_\pi)+2S_{\mathrm{nt}}(H_\pi)),
    \end{align}
    where
\begin{align*}
S_{\mathrm{tr}}(H_\pi)
&\triangleq \sum_{k\ge 2}\ \sum_{T\in \maT_k(H_\pi)} \inj(T,H_\pi)
\left(\frac{2e^{\kappa_b+1}}{n}\right)^k \big[c^{-1}(1+\gamma)\big]^{k-1},
\\
S_{\mathrm{nt}}(H_\pi)
&\triangleq \sum_{k\ge 2}\ \sum_{C\in \maC_k(H_\pi)} \inj(C,H_\pi)
\left(\frac{2e^{\kappa_b+1}}{n}\right)^k \big[c^{-1}(1+\gamma)\big]^{\xi k}.
\end{align*}
Here $\maT_k(H_\pi)$ and $\maC_k(H_\pi)$ denote the collections of connected trees and connected non-tree graphs with $k$ vertices that embed in $H_\pi$, respectively; $\inj(F,H_\pi)$ denotes the number of injective homomorphisms from $F$ to $H_\pi$ (see, e.g., \cite[Section 5.2]{lovasz2012large}). 
\end{proposition}

We first verify that the right-hand side of~\eqref{eq:lemma-str-snt} is $1+o(1)$. Note that $$\frac{r^2}{s^2}\le \frac{(1+\epsilon)^2 s^2}{n^2}\le \frac{(1+\epsilon)^2 (b\varrho^{-1}(b^{-1})-\delta)\log n}{In} = o(1),$$
    where the last equality follows from $p = n^{-a+o(1)},\gamma = n^{b+o(1)}, 0<b\le a$, and $2a-b<1$. Hence, $\tfrac{1}{1-e^{\kappa_b}r^2/s^2} = 1+o(1)$.
Following a similar argument as \cite[Lemma 3.5]{ding2023detection}, we have \begin{equation}\label{eq:counting-lemma}
\sum_{T\in \maT_k(H_\pi)} \inj(T,H_\pi)\le r(4\log n)^{k-1},
\quad
\sum_{C\in \maC_k(H_\pi)} \inj(C,H_\pi)\le k^3 n^{10\eta k}.
\end{equation}
Using \eqref{eq:counting-lemma} and $1+\gamma=n^{b+o(1)}$, we obtain
\begin{align}\label{eq:Str-o1}
S_{\mathrm{tr}}(H_\pi)
&\le \sum_{k\ge 2} r(4\log n)^{k-1}\Big(\frac{2e^{\kappa_b+1}}{n}\Big)^k \big[c^{-1}(1+\gamma)\big]^{k-1}\le \frac{Cr}{n}\sum_{k\ge 2}\big(C(\log n)n^{-(1-b)+o(1)}\big)^{k-1}=o(1),
\end{align}
because $b<1$ and $r\le n$.
Similarly,
\begin{align}\label{eq:Snt-o1}
S_{\mathrm{nt}}(H_\pi)
&\le \sum_{k\ge 2} k^3 n^{10\eta k}\Big(\frac{2e^{\kappa_b+1}}{n}\Big)^k \big[c^{-1}(1+\gamma)\big]^{\xi k}\le \sum_{k\ge 2} k^3 C^k n^{-(1-b\xi-10\eta+o(1))k}=o(1)
\end{align}
by $10\eta<1-b\xi$. We then prove Proposition~\ref{prop:main-IT-lwbd-poly-gamma}.

\subsection{Proof of Proposition~\ref{prop:prop-IT-lwbd-mid}}\label{subsec:proof-prop-IT-lwbd-mid}

We first establish the following lemma, which counts the number of permutations $\ti\pi$ for which $H(J_{\pi,\ti\pi})$ has a prescribed decomposition.
Let $F$ be a graph with no isolated vertices whose connected components consist of $x_i$ copies of distinct non-tree graphs $C_i$ and $y_j$ copies of distinct trees $T_j$. Let $v(F) = \sum_i x_i |C_i|+\sum_{j} y_j |T_j|$. For any $u\ge0$ and $\pi\in\Pi_r$, let $N_u(F,H_\pi)$ denote the
number of $\ti\pi\in\Pi_r$ such that the connected components of
$F$ consist of $x_i$ copies of $C_i$ and $y_j$
copies of $T_j$, and such that
$|S^*\setminus V(F)|=u$.

Fix $\ti\pi\in \Pi_r$ and denote $F = H(J_{\pi,\ti\pi})$. If the connected components of $F$ consist of $x_i$ copies of distinct non-tree graphs $C_i$ and $y_j$ copies of distinct trees $T_j$, then we have \begin{align}\label{eq:connected-comp-edge}
    |E(J_{\pi,\ti\pi})| = |E(F)|\le \sum_i x_i\xi |C_i|+\sum_j y_j(|T_j|-1),
\end{align}
because every connected non-tree subgraph of graph $H_\pi$ has edge-vertex ratio at most $\xi$ (see~\eqref{eq:ratio-xi}). For $u\ge 0$, let $F$ range over all possible active graphs
$H(J_{\pi,\ti\pi})$ with no isolated vertices. Since
$|S^*|=u+v(F)$ and $|E(J_{\pi,\ti\pi})|=|E(F)|$, we have
\begin{align*}
&\frac{1}{|\Pi_r|}
\sum_{\ti\pi\in\Pi_r}
e^{\kappa_b|S^*|}
[c^{-1}(1+\gamma)]^{|E(J_{\pi,\ti\pi})|} \quad=
\sum_{u\ge0}\sum_F
\frac{N_u(F;H_\pi)}{|\Pi_r|}
e^{\kappa_b(u+v(F))}
[c^{-1}(1+\gamma)]^{|E(F)|}.
\end{align*}
We then bound $N_u(F;H_\pi)$ by the following lemma.

\begin{lemma}\label{lem:decomp-HJ}
For every $u\ge 0$ and $\pi\in\Pi_r$, we have\begin{align*}
    N_u(F;H_\pi)\le |\Pi_r|\pth{\frac{r^2}{s^2}}^u\pth{\frac{2e}{n}}^{v(F)}\prod_{i} {[\inj(C_i,H_\pi)]^{x_i}}\prod_{j}{[\inj(T_j,H_\pi)]^{y_j}}
\end{align*}
    \end{lemma}
\begin{proof}

Let $\Aut(F)$ define the number of automorphisms of graph $F$. For each isomorphism class $C_i$, the number of ways to choose $x_i$ embedded copies in $H_\pi$ is at most
$(\inj(C_i,H_\pi)/\Aut(C_i))^{x_i}/x_i!$. Once these copies are chosen, the
restriction of $\theta=\ti\pi^{-1}\circ\pi$ may permute the $x_i$
copies and act by an automorphism within each copy, giving at most $x_i!\,\Aut(C_i)^{x_i}$
possible actions. Thus the total contribution of this isomorphism
class is at most 
$(\inj(C_i,H_\pi))^{x_i}.$
The same argument applies to the tree components. Therefore, the number of choices together with the action
of $\theta=\ti\pi^{-1}\circ\pi$ is bounded by
\begin{equation}\label{eq:embedding-choice}
\prod_i
\big[\inj(C_i,H_\pi)\big]^{x_i}
\prod_j
\big[\inj(T_j,H_\pi)\big]^{y_j}.
\end{equation}

Let $v = v(F)$. Fix $v$ distinct domain vertices $u_1,\cdots, u_v$ and $v$ distinct  range vertices $w_1,\cdots, w_v$, and prescribe $\ti\pi(u_i)  =w_i$ for $i = 1,\cdots, v$. Since $\ti\pi \sim \mathrm{Unif}(\Pi_r)$,\begin{align}\label{eq:prescribe-num}
    \prob{\ti\pi(u_i) = w_i,i=1,\cdots,v} = \frac{\binom{s-v}{r-v}^2 (r-v)!}{\binom{s}{r}^2 r!} = \frac{(r)_v}{(s)_v}\cdot \frac{1}{(s)_v}\le \Big(\frac{2e}{n}\Big)^v,
\end{align}
where the last inequality is because $\tfrac{(r)_v}{(s)_v}\le (\tfrac{r}{s})^v$, $\tfrac{1}{(s)_v}\le (\tfrac{e}{s})^v$, and $r \le 2 s^2/n$. After fixing $F$, the remaining part of $\ti\pi$ is a uniform partial injection of size $r-v$ between two sets of size $s-v$. A similar argument as~\eqref{eq:S^*-bound} in Lemma~\ref{lem:property_of_I} gives \begin{align*}
    \prob{|S^*\backslash V(J_{\pi,\ti\pi})| = u\mid V(J_{\pi,\ti\pi})} \le \pth{\frac{(r-v)^2}{(s-v)^2}}^{u}\le \pth{\frac{r^2}{s^2}}^u.
\end{align*}
Combining this with~\eqref{eq:embedding-choice} and~\eqref{eq:prescribe-num} yields Lemma~\ref{lem:decomp-HJ}.
\end{proof}
By Lemma~\ref{lem:decomp-HJ}, we have 
\begin{align}
\nonumber&~\frac{N_u(F;H_\pi)}{|\Pi_r|}
e^{\kappa_b(u+v(F))}
[c^{-1}(1+\gamma)]^{|E(F)|} \\\label{eq:sum-prop9}
\le&~
\left(e^{\kappa_b}\frac{r^2}{s^2}\right)^u
\left(\frac{2e^{\kappa_b+1}}{n}\right)^{v(F)}
[c^{-1}(1+\gamma)]^{|E(F)|} \times
\prod_i[\inj(C_i,H_\pi)]^{x_i}
\prod_j[\inj(T_j,H_\pi)]^{y_j}.
\end{align}
By~\eqref{eq:connected-comp-edge}, define
\begin{align*}
z_C
&\triangleq
\inj(C,H_\pi)
\left(\frac{2e^{\kappa_b+1}}{n}\right)^{|C|}
[c^{-1}(1+\gamma)]^{\xi|C|},\\
z_T
&\triangleq
\inj(T,H_\pi)
\left(\frac{2e^{\kappa_b+1}}{n}\right)^{|T|}
[c^{-1}(1+\gamma)]^{|T|-1}.
\end{align*}
It follows from~\eqref{eq:Str-o1} and~\eqref{eq:Snt-o1} that 
\begin{align}\label{eq:sum-zczt}
    \sum_T z_T=S_{\mathrm{tr}}(H_\pi)=o(1),\quad\sum_C z_C=S_{\mathrm{nt}}(H_\pi)=o(1).
\end{align}
The summand~\eqref{eq:sum-prop9} is therefore bounded by $\left(e^{\kappa_b}\frac{r^2}{s^2}\right)^u
\prod_i z_{C_i}^{x_i}\prod_j z_{T_j}^{y_j}.$
We now sum over $u$ and all possible component multiplicities.
Dropping the vertex-disjointness and total-size constraints only
increases the sum, and hence
\begin{align*}
&~\frac{1}{|\Pi_r|}
\sum_{\ti\pi\in\Pi_r}
e^{\kappa_b|S^*|}
[c^{-1}(1+\gamma)]^{|E(J_{\pi,\ti\pi})|}\\
\le&~
\sum_{u\ge0}
\left(e^{\kappa_b}\frac{r^2}{s^2}\right)^u
\prod_{k\ge2}\prod_{C\in\maC_k(H_\pi)}
\left(\sum_{x\ge0}z_C^x\right)
\prod_{k\ge2}\prod_{T\in\maT_k(H_\pi)}
\left(\sum_{y\ge0}z_T^y\right)\\
\le&~
\frac{1}{1-e^{\kappa_b}r^2/s^2}
\prod_{k\ge2}\prod_{C\in\maC_k(H_\pi)}
\frac{1}{1-z_C}
\prod_{k\ge2}\prod_{T\in\maT_k(H_\pi)}
\frac{1}{1-z_T},
\end{align*}
where the last inequality follows from $\sum_C z_C = o(1)$ and $\sum_{T}z_T = o(1)$.
Using $-\log(1-z)\le2z$ for $0\le z\le1/2$, we obtain
\[
\prod_C\frac{1}{1-z_C}\prod_T\frac{1}{1-z_T}
\le
\exp\!\Big(
2S_{\mathrm{nt}}(H_\pi)+2S_{\mathrm{tr}}(H_\pi)
\Big).
\]
This proves Proposition~\ref{prop:prop-IT-lwbd-mid}.

\section{Proof of Lemmas}
\subsection{Proof of Lemma~\ref{lem:orbit-expect}}\label{apd:proof-lem-orbit-expect}

For any $P = (e_1,\pi(e_1),e_2,\cdots,e_j,\pi(e_j))\in \sfP$ with $\ti{\pi}(e_2) = \pi(e_1),\cdots,\ti{\pi}(e_j) = \pi(e_{j-1})$, we have that \begin{align*}
    L_P &= \prod_{i=1}^j\ell(e_i,\pi(e_i)) \prod_{i=2}^{j}\ell(e_i,\ti{\pi}(e_i)) = \prod_{i=1}^j\ell(e_i,\pi(e_i)) \prod_{i=2}^{j}\ell(e_i,{\pi}(e_{i-1}))\\&=\ell(e_1,\pi(e_1))\ell(\pi(e_1),e_2)\cdots\ell(\pi(e_{j-1}),e_j)\ell(e_j,\pi(e_j)). 
\end{align*}
Therefore, $L_P$ can be expressed as $\ell(x,B_1)\ell(B_1,B_2)\cdots\ell(B_{k-1},B_k)\ell(B_k,y)$ for some $k\in \mathbb{N}$ and $x,B_1,\cdots,B_k,y\overset{\mathrm{i.i.d.}}{\sim} \Bernoulli(p)$. Denote $$\ell^{(k)}_\sfP(x,y)\triangleq \mathbb{E}_{\maQ}\qth{\ell(x,B_1)\pth{\prod_{i=1}^{k-1}\ell(B_i,B_{i+1})}\ell(B_k,y)},$$ then $\mathbb{E}_{\maQ}[L_P] = \mathbb{E}_{{\maQ}}\qth{\ell^{(k)}(x,y)}$. We note that 
     \begin{align*}
        \ell^{(k)}_\sfP(x,y) &= \mathbb{E}_{\maQ}\qth{\ell(x,B_1)\pth{\prod_{i=1}^{k-1}\ell(B_i,B_{i+1})}\ell(B_k,y)}\\
        &=\sum_{b_1,\cdots,b_{k}\in \{0,1\}} \ell(x,b_1)\prob{B_1= b_1}\pth{\prod_{i=1}^{k-1} \ell(b_i,b_{i+1}) \prob{B_{i+1} = b_{i+1}}}\ell(b_{k},y).
    \end{align*}
Define $M(a,b)\triangleq \ell(a,b)\prob{B=b}$ for $B\sim \mathrm{Bern}(p)$ and $a,b\in \{0,1\}$ and a matrix \begin{align*}
    M\triangleq \begin{bmatrix}
        M(0,0) & M(0,1) \\ M(1,0) & M(1,1)
    \end{bmatrix} = \begin{bmatrix}
        1-p+p\rho & (1-p)(1-\rho)\\ p(1-\rho) & p+\rho-p\rho
    \end{bmatrix}.
\end{align*}
Then, we obtain that \begin{align*}
    \ell^{(k)}_\sfP(x,y) &= \sum_{b_1,\cdots,b_{k}\in \{0,1\}} M(x,b_1)M(b_1,b_2)\cdots M(b_{k-1},b_{k})\ell(b_{k},y) \\&= \sum_{b_{k}\in \{0,1\}} M^{k}(x,b_{k}) \ell(b_{k},y),
\end{align*}
which is of form $\ell^{(k)}_\sfP(x,y) = \alpha(x,y) \lambda_1^{k}+\beta(x,y)\lambda_2^k$, where $\lambda_1,\lambda_2$ is the eigenvalue of matrix $M$. By direct calculation, we obtain that $\lambda_1 = 1,\lambda_2 = \rho$. The coefficients $\alpha(x,y)$ and $\beta(x,y)$ can be determined via the first two terms $k=1,2$. Then we get \begin{align*}
    \ell^{(k)}_\sfP(1,1) = 1+\frac{1-p}{p}\cdot \rho^{k+1},\ell^{(k)}_\sfP(0,1) = \ell^{(k)}_\sfP(1,0) = 1-\rho^{k+1},\ell^{(k)}_\sfP(0,0) = 1+\frac{p}{1-p}\cdot \rho^{k+1}.
\end{align*}
Therefore, we obtain that \begin{align*}
    \mathbb{E}_{\maQ}\qth{\ell^{(k)}_\sfP(x,y)} &= p^2 \ell^{(k)}_\sfP(1,1)+p(1-p)\ell^{(k)}_\sfP(0,1)\\&~~~+p(1-p)\ell^{(k)}_\sfP(1,0)+(1-p)^2\ell^{(k)}_\sfP(0,0)=1,
\end{align*}
yielding that $\mathbb{E}_{\maQ}[L_P] = 1$ for any $P\in \sfP^\sfE$.

For any $C = (e_1,\pi(e_1),e_2,\cdots,e_j,\pi(e_j))\in \sfC^\sfE$ with $\ti{\pi}(e_2) = \pi(e_1),\cdots,\ti{\pi}(e_j) = \pi(e_{j-1})$ and $\ti{\pi}(e_1) = \pi(e_j)$, denote $e_0 =e_j$, then we have \begin{align*}
    L_C &= \prod_{i=1}^j\ell(e_i,\pi(e_i)) \prod_{i=1}^{j}\ell(e_i,\ti{\pi}(e_i)) = \prod_{i=1}^j\ell(e_i,\pi(e_i)) \prod_{i=1}^{j}\ell(e_i,{\pi}(e_{i-1}))\\&=\ell(e_1,\pi(e_1))\ell(\pi(e_1),e_2)\cdots\ell(\pi(e_{j-1}),e_j)\ell(e_j,\pi(e_j))\ell(\pi(e_j),e_1). 
\end{align*}
Therefore, $L_C$ can be expressed as $\ell(B_1,B_2)\cdots\ell(B_{k-1},B_k)\ell(B_k,B_1)$ for  $k=2j$ and $B_1,\cdots,B_k\overset{\mathrm{i.i.d.}}{\sim} \Bernoulli(p)$. 
Denote $B_{k+1} = B_1$, we note that 
     \begin{align*}
        \mathbb{E}_{\maQ}[L_C] &= \mathbb{E}_{\maQ}{\qth{\prod_{i=1}^{k}\ell(B_i,B_{i+1})}}=\sum_{b_1,\cdots,b_{k}\in \{0,1\}} {\prod_{i=1}^{k} \ell(b_i,b_{i+1}) \prob{B_{i+1} = b_{i+1}}}\\&=\sum_{b_1,\cdots,b_{k}\in \{0,1\}} \prod_{i=1}^{k} M(b_i,b_{i+1}) = \sum_{b_1\in\{0,1\}} M^{k}(b_1,b_1) = \mathrm{tr}(M^k).
    \end{align*}
    Since $\mathrm{tr}(M^k) = \lambda_1^k+\lambda_2^k = 1+\rho^{k}$ and $k = 2|C|$, we obtain that $\mathbb{E}_{\maQ}[L_C] = 1+\rho^{2|C|}$.

\subsection{Proof of Lemma~\ref{lem:property_of_I}}\label{apd:proof-lem-property-I}
Since $\pi$ and $\ti\pi$ are uniformly distributed over injective mappings satisfying
$|\sfS_\pi|=|\sfT_\pi|=r$, by symmetry we may fix $\pi$ and take the randomness only over $\ti\pi$. 
Let $ K\triangleq\sum_{k\ge1}k a_k$. If \(K>r\), then no such collection of vertex-disjoint cycles can occur and
the desired bound is trivial. Hence we assume \(K\le r\). For each \(k\), the falling factorial \((C_k)_{a_k}\) counts ordered \(a_k\)-tuples of distinct \(k\)-cycles appearing in \(\theta\triangleq \ti\pi^{-1}\circ \pi\) (see Figure~\ref{fig:vertex-cycles} for an illustration). 
 For a cyclic ordering $c=(v_1,\ldots,v_k)$ of $k$ distinct vertices, say
that $c$ is realized if
\[
\theta(v_i)=v_{i+1}\quad (1\le i\le k),
\qquad v_{k+1}=v_1.
\]
Let $\mathfrak C$ range over ordered collections containing $a_k$ distinct
candidate $k$-cycles for each $k$. Then
\[
\expect{\prod_{k\ge1}(C_k)_{a_k}}
=
\sum_{\mathfrak C}
\mathbb P\{\mathfrak C\text{ is realized}\}.
\]
Since $\theta$ is a partial injection, its directed cycles are
vertex-disjoint. Hence any collection $\mathfrak C$ containing two cycles
with a common vertex cannot be realized and contributes zero to the sum.
We may therefore restrict the sum to vertex-disjoint collections.

\begin{figure}[htbp]
\centering
\begin{tikzpicture}[
    x=1cm,y=1cm,
    dot/.style={circle,fill=black,minimum size=3.2pt,inner sep=0pt},
    vertex/.style={circle,draw=black!70,fill=white,line width=0.7pt,
        minimum size=6.8mm,inner sep=0pt,font=\small},
    pimap/.style={-{Stealth[length=2mm,width=1.35mm]},
        line width=0.75pt,draw=blue!70!black,shorten <=2.2pt,shorten >=2.2pt},
    tpimap/.style={-{Stealth[length=2mm,width=1.35mm]},
        line width=0.75pt,dashed,draw=red!70!black,
        shorten <=2.2pt,shorten >=2.2pt},
    reduce/.style={-{Stealth[length=2mm,width=1.35mm]},
        line width=0.7pt,draw=black!75},
    orbit/.style={-{Stealth[length=2mm,width=1.4mm]},
        line width=0.85pt,draw=black!80},
    note/.style={font=\small,text=black!80}
]
\node[font=\scriptsize,anchor=east] at (0.25,1.15)
    {$\sfS_\pi\cap\sfS_{\ti\pi}$};
\node[font=\scriptsize,anchor=east] at (0.25,-0.15)
    {$\sfT_\pi\cap\sfT_{\ti\pi}$};

\node[dot,label=above:$v_1$] (s1) at (1.20,1.15) {};
\node[dot,label=above:$v_2$] (s2) at (3.20,1.15) {};
\node[dot,label=below:$\pi(v_1)$] (t1) at (1.20,-0.15) {};
\node[dot,label=below:$\pi(v_2)$] (t2) at (3.20,-0.15) {};
\draw[pimap] (s1)--(t1)
    node[left,font=\scriptsize,pos=0.43] {$\pi$};
\draw[pimap] (s2)--(t2);
\draw[tpimap] (s2)--(t1)
    node[above,font=\scriptsize,pos=0.50] {$\ti\pi$};
\draw[tpimap] (s1)--(t2);
\draw[reduce] (2.20,-0.72)--(2.20,-1.2);
\node[vertex] (c21) at (1.30,-2) {$v_1$};
\node[vertex] (c22) at (3.10,-2) {$v_2$};
\draw[orbit] (c21) to[bend left=30] (c22);
\draw[orbit] (c22) to[bend left=30] (c21);
\node[note] at (2.20,-2.95) {$C_2=1$};

\node[dot,label=above:$v_3$] (s3) at (6.20,1.15) {};
\node[dot,label=above:$v_4$] (s4) at (8.20,1.15) {};
\node[dot,label=above:$v_5$] (s5) at (10.20,1.15) {};
\node[dot,label=below:$\pi(v_3)$] (t3) at (6.20,-0.15) {};
\node[dot,label=below:$\pi(v_4)$] (t4) at (8.20,-0.15) {};
\node[dot,label=below:$\pi(v_5)$] (t5) at (10.20,-0.15) {};
\draw[pimap] (s3)--(t3)node[left,font=\scriptsize,pos=0.43] {$\pi$};
\draw[pimap] (s4)--(t4);
\draw[pimap] (s5)--(t5);
\draw[tpimap] (s4)--(t3);
\draw[tpimap] (s5)--(t4)
    node[above,font=\scriptsize,pos=0.43] {$\ti\pi$};
\draw[tpimap] (s3)--(t5);
\draw[reduce] (8.20,-0.72)--(8.20,-1.2);
\node[vertex] (c33) at (6.85,-2.40) {$v_3$};
\node[vertex] (c34) at (8.20,-1.62) {$v_4$};
\node[vertex] (c35) at (9.55,-2.40) {$v_5$};
\draw[orbit] (c33)--(c34);
\draw[orbit] (c34)--(c35);
\draw[orbit] (c35)--(c33);
\node[note] at (8.20,-2.95) {$C_3=1$};
\end{tikzpicture}
\caption{Examples of vertex cycles of $\ti\pi^{-1}\circ\pi$. Solid and
dashed arrows denote $\pi$ and $\ti\pi$, respectively; the resulting
components contribute one to $C_2$ and $C_3$, respectively.}
\label{fig:vertex-cycles}
\end{figure}

We now count vertex-disjoint ordered candidate collections. Order the cycles first by their length \(k\), and within each fixed \(k\) by the order appearing in the falling factorial. Choose an ordered list of \(K\) distinct vertices from \(\sfS_\pi\), which gives \((r)_K\) possibilities. Split this ordered list into successive blocks of lengths
\[
    \underbrace{1,\ldots,1}_{a_1\text{ times}},
    \underbrace{2,\ldots,2}_{a_2\text{ times}},
    \underbrace{3,\ldots,3}_{a_3\text{ times}},
    \ldots .
\]
A block \((v_1,\ldots,v_k)\) represents the directed cycle
\[
    v_1\mapsto v_2,\quad v_2\mapsto v_3,\quad \ldots,\quad
    v_k\mapsto v_1
\]
such that $\pi(v_1) = \ti{\pi}(v_2),\pi(v_2) = \ti{\pi}(v_3),\cdots,\pi(v_k) = \ti{\pi}(v_1)$.
The same directed \(k\)-cycle is represented by exactly \(k\) cyclic rotations
of the block. Hence the number of vertex-disjoint ordered candidate
collections is exactly $\tfrac{(r)_K}{\prod_{k\ge1}k^{a_k}}.$

Fix one such candidate collection \(\mathfrak C\).
Without loss of generality we can fix $\pi$ first.
Realizing it imposes
exactly \(K\) prescribed mapping constraints.
where the \(u_i\)'s are distinct and the \(v_i\)'s are distinct. The total
number of size-\(r\) partial injections from \(\sfS_\pi\) to \(\sfT_\pi\) is $M_r=\binom{s}{r}^2 r!$.
If the \(K\) prescribed pairs are imposed, the remaining \(r-K\) domain
vertices and \(r-K\) range vertices can be chosen in $\binom{s-K}{r-K}^2$
ways and then bijected in \((r-K)!\) ways. Hence
\begin{align*}
    \mathbb P\qth{\mathfrak C\text{ is realized}}
    &=
    \frac{\binom{s-K}{r-K}^2 (r-K)!}{\binom{s}{r}^2 r!} =
    \frac{(r)_K}{(s)_K^2}.
\end{align*}
Therefore
\begin{align*}
    \expect{\prod_{k\ge1}(C_k)_{a_k}}
    &\le
    \frac{(r)_K}{\prod_{k\ge1}k^{a_k}}
    \cdot
    \frac{(r)_K}{(s)_K^2} =
    \frac{(r)_K^2}{(s)_K^2}
    \prod_{k\ge1} k^{-a_k}.
\end{align*}
Since \((r)_K/(s)_K\le (r/s)^K\), we obtain
\begin{align*}
    \expect{\prod_{k\ge1}(C_k)_{a_k}}
    \le
    \left(\frac{r}{s}\right)^{2K}
    \prod_{k\ge1} k^{-a_k} = \prod_{k\ge 1}\pth{\frac{(r^2/s^2)^k}{k}}^{a_k}.
\end{align*}

We then derive the fixed-point bound. Applying~\eqref{eq:cycle-moment} with \(a_1=f\) and \(a_k=0\) for \(k\ge2\) gives $\expect{(C_1)_{c_1}}\le (r^2/s^2)^{c_1}.$
For \(c_1\ge1\), since $\indc{C_1\ge c_1}\le \binom{C_1}{c_1} = \tfrac{(C_1)_{c_1}}{c_1!}$, we obtain 
\begin{align*}
    \prob{C_1 = c_1}\le \prob{C_1\ge c_1}\le \frac{\expect{(C_1)_{c_1}}}{c_1!}\le \frac{(r^2/s^2)^{c_1}}{c_1!}
\end{align*}
Similarly, using $\indc{C_k = a_k}\le \tbinom{C_k}{a_k}\le \tfrac{(C_k)_{a_k}}{a_k!}$, we have \begin{align*}
    \prob{|S^*(\pi,\ti{\pi})| = t} = 
    \prob{\sum_{k\ge 1}kC_k = t}\le \sum_{\sum_{k\ge 1} ka_k = t} \prod_{k\ge 1} \frac{1}{a_k!}\pth{\frac{(r^2/s^2)^{k}}{k}}^{a_k} = (r^2/s^2)^t,
\end{align*}
where the last equality follows by extracting the coefficient of $z^t$ in the exponential generating function
\[
\prod_{k\ge1}\exp\Big(\frac{(r^2/s^2)^k z^k}{k}\Big)
=
\exp\Big(\sum_{k\ge1}\frac{((r^2/s^2)z)^k}{k}\Big)
= \frac{1}{1-(r^2/s^2)z}.
\]
Indeed, the coefficient of $z^t$ in the last display is exactly $(r^2/s^2)^t$. This proves the lemma.

\subsection{Proof of Lemma~\ref{lem:conditional-prob-lwbd}}\label{subsec:proof-lem-condition-prob}

Note that $H(J)$ is a subgraph of $H_\pi$. Since $H_\pi\in \maH(r)$, we have $H(J)\in \maH(|V(J)|)$.
Conditional on $E(J)\subseteq E(H_\pi)$, all remaining edge indicators
of $H_\pi$ are independent $\operatorname{Bern}(p_{11})$ variables. Moreover,
\[
\{J_{\pi,\ti\pi}=J\}
=
\{E(J)\subseteq E(H_\pi)\}\cap\mathcal D_J,
\]
where
\[
\mathcal D_J
\triangleq
\bigcap_{C\in\sfC^\sfE\setminus J}
\left\{
C\cap\binom{\sfS_\pi}{2}
\not\subseteq E(H_\pi)
\right\}
\]
is decreasing in the remaining edges. The first event ensures that all edges of $J$ are present in $H_\pi$,
whereas $\mathcal D_J$ requires every cycle
$C\in\sfC^\sfE\setminus J$ to contain at least one missing edge,
thereby excluding any additional cycle.

Recall that $\eta$ is a constant such that $10\eta<1-b\xi$.
Let $K_0=\lceil4/\eta\rceil$. Conditional on
$E(J)\subseteq E(H_\pi)$, let $A_1$ be the event that no edge outside
$E(J)$ is incident to $V(J)$, and let $A_2$ be the event that $H_\pi[\sfS_\pi\setminus V(J)]$
satisfies the density and maximum-degree conditions in the definition
of $\maH(r)$ and contains no cycle of length $3\le k<K_0$. Then both $A_1$ and $A_2$ are decreasing events.
Recall that $rp_{11}\le \xi$. Since $A_1$ requires
at most $|V(J)|r$ edges to be absent,
\[
\maP_r(A_1\mid\pi,E(J)\subseteq E(H_\pi))
\ge (1-p_{11})^{|V(J)|r}
\ge e^{-2\xi|V(J)|}.
\]
The density and maximum-degree conditions in $A_2$ hold with probability
$1-o(1)$ since $\maP_r(\maE_\pi^2) = 1-o(1)$. In addition, FKG inequality gives
\[
\maP_r\!\left(
H_\pi[U]\text{ has no cycle of length }3\le k<K_0
\right)
\ge
\prod_{k=3}^{K_0-1}(1-p_{11}^k)^{r^k}
\ge a_0,
\]
where
\[
a_0\triangleq
\exp\Big(-2\sum_{k=3}^{K_0-1}\xi^k\Big)>0.
\]
Hence another application of FKG inequality yields
\[
\maP_r(A_2\mid\pi,E(J)\subseteq E(H_\pi))
\ge (1-o(1))a_0.
\]

On $A_1\cap A_2$, the graph $H_\pi$ is the disjoint union of $H(J)$
and $H_\pi[\sfS_\pi\setminus V(J)]$. The density and maximum-degree conditions therefore
hold. For $3\le k<K_0$, all $k$-cycles lie in the graph
$H(J)$, while for $k\ge K_0$,
\[
\#\{\text{labeled $k$-cycles in }H_\pi\}
\le r(\log n)^k
\le n^{1+\eta k/4}
\le n^{\eta k}.
\]
Thus $A_1\cap A_2\subseteq\maE_\pi^2$.
Finally, $A_1\cap A_2$ and $\mathcal D_J$ are decreasing. By FKG and
the independence of the edge sets involved in $A_1$ and $A_2$,
\begin{align*}
\maP_r(\maE_\pi^2\mid\pi,J_{\pi,\ti\pi}=J)
&\ge
\maP_r(A_1\cap A_2
 \mid\pi,E(J)\subseteq E(H_\pi),\mathcal D_J)\ge
(1-o(1))a_0e^{-2\Lambda|V(J)|}.
\end{align*}
If $J\ne\emptyset$, then $|V(J)|\le2|E(J)|$, so after decreasing the
constant if necessary, the last display is at least $c^{|E(J)|}$ for some
$c>0$. If $J=\emptyset$, FKG gives
\[
\maP_r(\maE_\pi^2\mid\pi,J_{\pi,\ti\pi}=\emptyset)
\ge \maP_r(\maE_\pi^2\mid\pi)=1-o(1).
\]
Consequently,
\[
\maP_r(\maE_\pi^2\mid\pi,J_{\pi,\ti\pi}=J)
\ge (1-o(1))c^{|E(J)|},
\]
which proves the lemma.

\section{Auxiliary Lemmas}

\begin{lemma}\label{lem:var-gamma-bounded}
    For any $\gamma<\infty$ and $e\in \binom{\sfS_\pi}{2}$, we have \begin{align*}
        \var_{\maP_{r}}[\log\ell(G_1(e),G_2(\pi(e)))]\le (2+\log(1+\gamma))I\le (2+\gamma)I.
    \end{align*}
\end{lemma}
\begin{proof}
Let $R=\ell(G_1(e),G_2(\pi(e)))$. Note that
\[
\var_{\maP_{r}}(\log R)\le \mathbb E_{\maP_{r}}[(\log R)^2]=\mathbb E_{\maQ}[R(\log R)^2].
\]
For $0\le x\le 1+\gamma$, we claim that
\[
x(\log x)^2\le (2+\gamma)(x\log x-x+1).
\]
Indeed, for $0\le x\le1$, one has $x(\log x)^2\le 2(x\log x-x+1)$. For $0\le x\le1$, write $x=e^{-t}$ with $t\ge0$. Then
\[
x(\log x)^2= e^{-t}t^2,\qquad
2(x\log x-x+1)=2(1-e^{-t}t-e^{-t}).
\]
Thus the desired inequality is equivalent to
\[
e^{-t}(t^2+2t+2)\le2,
\]
which follows because its left-hand side is decreasing in $t\ge0$ and equals $2$ at $t=0$.

For $1\le x\le1+\gamma$, writing $x=e^t$ gives
\[
\frac{x(\log x)^2}{x\log x-x+1}=\frac{t^2}{t-1+e^{-t}}.
\]
The right-hand side is increasing for $t\ge0$, since its derivative has the same sign as $t-2+(t+2)e^{-t}\ge0$. Therefore,
\[
\frac{x(\log x)^2}{x\log x-x+1}\le \frac{\log^2(1+\gamma)}{\log(1+\gamma)-1+(1+\gamma)^{-1}}\overset{\mathrm{(a)}}{\le} \log(1+\gamma)+2\le \gamma+2.
\]
where $\mathrm{(a)}$ follows from the fact that, with $L=\log(1+\gamma)$,
\[
L^2\le (L+2)(L-1+e^{-L}),
\]
which is equivalent to $L-2+(L+2)e^{-L}\ge0$; the latter holds since its derivative is
$1-(L+1)e^{-L}\ge0$ and it vanishes at $L=0$.

Taking expectation under $\maQ$ yields
\[
\mathbb E_{\maQ}[R(\log R)^2]\le (2+\gamma)\mathbb E_{\maQ}[R\log R-R+1]=(2+\gamma)I,
\]
where the last equality follows from $\mathbb E_{\maQ}R=1$ and $R=d\maP_r/d\maQ$. Combining the above inequalities proves the lemma.
\end{proof}

\begin{lemma}[Bernstein inequality]\label{lem:one-sided-bernstein}
Let $\xi_1,\ldots,\xi_N$ be independent random variables such that
$\mathbb E \xi_i=0$ and $\xi_i\le B$ almost surely for every $i$, where
$B>0$. Let
    $V=\sum_{i=1}^N \mathbb E \xi_i^2.$
Then, for every $x>0$,
\[
    \mathbb P\left[\sum_{i=1}^N \xi_i\ge x\right]
    \le
    \exp\Bigg(
        -\frac{x^2}{2(V+Bx/3)}
    \Bigg).
\]
\end{lemma}

\begin{proof}
The proof follows from~\cite[Chapter 2.8]{boucheron2013concentration}.
\end{proof}

\begin{lemma}\label{lem:psi}
Let $u_k\ge 0$ be finitely supported for $k\ge 2$, and define
\[
\Psi(u)=\sum_{k\ge2}\frac{(r^2/s^2)^k}{k}\left(e^{u_k}-1\right).
\]
Then, for every integer $c_1\ge0$,
\begin{align}\label{eq:lem-Psi-1}
    \mathbb E \left[
\indc{C_1=c_1}\exp\Big(\sum_{k\ge2}u_k C_k\Big)
\right]
\le
\frac{(r^2/s^2)^{c_1}}{c_1!}e^{\Psi(u)},
\end{align}
and
\begin{align}\label{eq:lem-Psi-2}
    \mathbb E \left[
\indc{C_1=c_1}\left(\exp\Big(\sum_{k\ge2}u_k C_k\Big)-1\right)
\right]
\le
\frac{(r^2/s^2)^{c_1}}{c_1!}\left(e^{\Psi(u)}-1\right).
\end{align}
The same bounds hold for countably many $u_k$ whenever the right-hand side is finite.
\end{lemma}

\begin{proof}
For each $k$,
\[
e^{u_k C_k}
=
\left(1+(e^{u_k}-1)\right)^{C_k}
=
\sum_{a_k\ge0}\frac{(C_k)_{a_k}}{a_k!}(e^{u_k}-1)^{a_k}.
\]
Multiplying the above identity over $k\ge2$ and using monotone convergence reduces the proof to finitely many $k$'s. Also,
\[
\indc{C_1=c_1}\le \frac{(C_1)_{c_1}}{{c_1}!}.
\]
Therefore, by Lemma~\ref{lem:property_of_I}, term by term,
\[
\E\left[
\indc{C_1=c_1}\prod_{k\ge2}(C_k)_{a_k}
\right]
\le
\frac{1}{c_1!}\E\left[
(C_1)_{c_1}\prod_{k\ge2}(C_k)_{a_k}
\right]
\le
\frac{(r^2/s^2)^{c_1}}{c_1!}\prod_{k\ge2}\left(\frac{(r^2/s^2)^k}{k}\right)^{a_k}.
\]
Substitution into the expansion gives
\[
\begin{aligned}
\E\left[
\indc{C_1=c_1}\exp\Big(\sum_{k\ge2}u_k C_k\Big)
\right]
&\le
\frac{(r^2/s^2)^{c_1}}{c_1!}
\prod_{k\ge2}\sum_{a_k\ge0}
\frac{1}{a_k!}
\Big(\frac{(r^2/s^2)^k}{k}(e^{u_k}-1)\Big)^{a_k}  \\
&=
\frac{(r^2/s^2)^{c_1}}{c_1!}
\exp\Bigg(
\sum_{k\ge2}\frac{(r^2/s^2)^k}{k}(e^{u_k}-1)
\Bigg)
=
\frac{(r^2/s^2)^{c_1}}{c_1!}e^{\Psi(u)}.
\end{aligned}
\]
If the empty multi-index $(a_k)\equiv0$ is removed from the expansion, the same calculation gives
\[
\E\left[
\indc{C_1=c_1}\left(\exp\Big(\sum_{k\ge2}u_k C_k\Big)-1\right)
\right]
\le
\frac{(r^2/s^2)^{c_1}}{c_1!}\left(e^{\Psi(u)}-1\right).
\]
This proves the lemma.
\end{proof}

\bibliographystyle{alpha}
\bibliography{main}
\end{document}